\documentclass[final,3p,times]{elsarticle}

\usepackage{amssymb}
\usepackage{amsmath}
 \usepackage{amsthm, xspace}

\usepackage[T1]{fontenc}
\usepackage{tikz, mdframed, subcaption, amsmath, multirow, tikz, graphicx, lmodern, enumitem}
\usetikzlibrary{tikzmark,calc}

\newtheorem{theorem}{Theorem}[section]
\newtheorem{lemma}[theorem]{Lemma}
\newtheorem{corollary}[theorem]{Corollary}

\theoremstyle{definition}
\newtheorem{reduction}{Reduction}

\newcommand{\codes}{{\textsc{Codes}}}

\newcommand{\dd}{\gamma^{\rm{D}}}
\newcommand{\td}{\gamma^{\rm{TD}}}

\newcommand{\sset}{\gamma^{\rm{S}}}
\newcommand{\lset}{\gamma^{\rm{L}}}
\newcommand{\oset}{\gamma^{\rm{O}}}
\newcommand{\iset}{\gamma^{\rm{I}}}
\newcommand{\fset}{\gamma^{\rm{F}}}

\newcommand{\sd}{\gamma^{\rm{SD}}}
\newcommand{\std}{\gamma^{\rm{STD}}}

\newcommand{\x}{\gamma^{\rm{X}}}
\newcommand{\id}{\gamma^{\rm{ID}}}
\newcommand{\itd}{\gamma^{\rm{ITD}}}
\newcommand{\fd}{\gamma^{\rm{FD}}}
\newcommand{\ftd}{\gamma^{\rm{FTD}}}
\newcommand{\ld}{\gamma^{\rm{LD}}}

\newcommand{\od}{\gamma^{\rm{OD}}}

\newcommand{\hyp}{\mathcal{H}}
\newcommand{\hedge}{\mathcal{F}}
\newcommand{\clu}{\mathcal{C}}

\newcommand{\calO}{\mathcal{O}\xspace}
\newcommand{\calT}{\mathcal{T}\xspace}

\newcommand{\minset}[1]{\textsc{Min #1-Set}\xspace}

\newcommand{\test}{\textsc{Test Cover}\xspace}

\newcommand{\sadmis}{{\textsc{S}}-admissible}

\newcommand{\defproblem}[3]{
\vspace{3mm}
\noindent \fbox{
\begin{minipage}{0.96\textwidth}
\begin{tabular*}{\textwidth}{@{\extracolsep{\fill}}lr}
#1 \\ \end{tabular*}
{\bf{Input:}} #2 \\
{\bf{Question:}} #3
\end{minipage}
}\vspace{3mm}
}

\newcommand{\dc}[1]{{\color{blue}#1}}
\newcommand{\aw}[1]{{\color{red}#1}}

\usepackage{tikz}

\usepackage{pgffor} 

\usetikzlibrary{shapes.geometric}

\usepackage{intcalc}
\usetikzlibrary{automata, positioning}

\usetikzlibrary{calc}

\usepackage[thinlines]{easytable}

\usepackage[figuresright]{rotating}

\begin{document}

\begin{frontmatter}




\title{The Interplay Between Domination and Separation in Graphs}


 \author[label1,label2]{Dipayan Chakraborty}
 \author[label1]{Annegret K. Wagler}
 
 \address[label1]{LIMOS, Universit\'{e} Clermont Auvergne, CNRS, Mines Saint-\'{E}tienne, Clermont Auvergne INP, LIMOS, Clermont-Ferrand, France.}
 
 \address[label2]{Department of Computer Science and Mathematics, Lebanese American University, Beirut, Lebanon
 }



\begin{abstract}
In the literature, several identification problems in graphs have been studied, of which, the most widely studied are the ones based on dominating sets as a tool of identification. Hereby, the objective is to separate any two vertices of a graph by their unique neighborhoods in a suitably chosen dominating or total-dominating set. Such a (total-)dominating set endowed with a separation property is often referred to as a code of the graph. In this paper, we study the four separation properties location, closed-separation, open-separation and full-separation. We address the complexity of finding minimum separating sets in a graph and study the interplay of these separation properties with several codes (establishing a particularly close relation between separation and codes based on domination) as well as the interplay of separation and complementation (showing that location and full-separation are the same on a graph and its complement, whereas closed-separation in a graph corresponds to open-separation in its complement).
\end{abstract}


\end{frontmatter}


\section{Introduction}

\begin{figure}[!t]
\begin{center}
\includegraphics[scale=1.0]{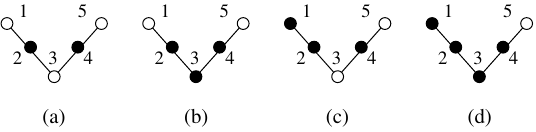}
\caption{Minimum S-sets and X-codes in a graph (the black vertices belong to the S-sets and X-codes), where (a) is both a D-set and an L-set and thus an LD-code, (b) is both a TD-set and an I-set, hence also both an ID- and an ITD-code, (c) is an O-set and thus an OD-code, (d) is an F-set, hence both an FD- and an FTD-code as well as an ITD-code.}
\label{fig_exp_X-codes}
\end{center}
\end{figure}

In the area of identification problems, the aim is to distinguish any two vertices of a graph by the unique intersection of their neighborhoods with a suitably chosen dominating or total-dominating set of the graph. Consider a simple, undirected graph $G=(V,E)$ and denote by $N(v)$ the \emph{open neighborhood} of a vertex $v \in V$. Moreover, let $N[v] = N(v) \cup \{v\}$ denote the \emph{closed neighborhood} of $v$. A subset $C \subseteq V$ is a
\begin{itemize}
  \itemsep-3pt
\item \emph{dominating set} (for short \emph{D-set}) if $N[v]\cap C$ are non-empty sets for all $v \in V$,
\item  \emph{total-dominating set} (for short \emph{TD-set}) if $N(v)\cap C$ are non-empty sets for all $v \in V$. 
\end{itemize}
To distinguish vertices of a graph, 
different separation properties 
have been studied. A subset $C \subseteq V$ is a
\begin{itemize}
  \itemsep-3pt
\item \emph{locating set} 
  (for short \emph{L-set}) if $N(v)\cap C$ is a unique set for each $v \in V\setminus C$, 
\item \emph{open-separating set} (for short \emph{O-set}) if $N(v)\cap C$ is a unique set for each $v \in V$,
\item \emph{closed-separating set} (for short \emph{I-set}) if $N[v]\cap C$ is a unique set for each $v \in V$, and
\item \emph{full-separating set} (for short \emph{F-set}) if it is both open- and closed-separating, that is if $N[u] \cap C \ne N[v] \cap C$ and $N(u) \cap C \ne N(v) \cap C$ holds for each pair of distinct vertices $u,v \in V$.
\end{itemize}
For illustration, Figure \ref{fig_exp_X-codes} exhibits in a small graph examples of such sets of minimum cardinality.

Combining one separation property with one domination property, the following identification problems have been studied in the literature:
\begin{itemize}
  \itemsep-3pt
\item location with domination and total-domination leading to \emph{locating dominating codes} (for short \emph{LD-codes}) and \emph{locating total-dominating codes} (for short \emph{LTD-codes}), see~\cite{S_1988} and \cite{HHH_2006}, respectively;
\item open-separation with domination and total-domination leading to \emph{open-sepa\-ra\-ting dominating codes} (\emph{OD-codes} for short) and \emph{open-separating total-do\-minating codes} \footnote{OTD-codes were introduced independently in~\cite{HLR_2002} and in~\cite{SS_2010} under the names of \emph{strongly $(t,\le l)$-identifying codes} and \emph{open neighborhood locating-dominating sets}   (or \emph{OLD-sets}), respectively. However, due to consistency in naming that specifies the separation and the domination property, we prefer to call them open-separating total-dominating codes in this article.}
  (\emph{OTD-codes} for short), see \cite{CW_ISCO2024} and  \cite{SS_2010}, respectively;
\item closed-separation with domination and total-domination leading to \emph{identifying codes} (\emph{ID-codes} for short) and \emph{identifying total-dominating codes}\footnote{Identifying total-dominating codes had been introduced to the literature in~\cite{HHH_2006} under the name \emph{differentiating total-dominating codes}. However, due to consistency in notation, we prefer to call them ITD-codes in this article.} 
  (\emph{ITD-codes} for short), see \cite{KCL_1998} and \cite{HHH_2006}, respectively;
\item full-separation with domination and total-domination leading to \emph{full-sepa\-ra\-ting dominating codes} (for short \emph{FD-codes}) and \emph{full-separating total-do\-minating codes} (for short \emph{FTD-codes}), respectively, see \cite{CW_Caldam2025}.
\end{itemize}
Figure \ref{fig_exp_X-codes} illustrates examples of such codes in a small graph. 
It is known from the literature that not every graph has all studied types of sets and codes. 
Summarizing results from for instance \cite{CW_ISCO2024,CW_Caldam2025,HHH_2006,KCL_1998,SS_2010},
we see that a graph $G$ has 
\begin{itemize}
  \itemsep-3pt
\item no TD-set if $G$ has \emph{isolated vertices}, that is, vertices $v$ with $N(v) = \emptyset$,
\end{itemize}
whereas every graph has a D-set; a graph $G$ has 
\begin{itemize}
  \itemsep-3pt
\item no O-set if $G$ has \emph{open twins}, that is, non-adjacent vertices $u,v$ with $N(u) = N(v)$,
\item no I-set if $G$ has \emph{closed twins}, that is, adjacent vertices $u,v$ with $N[u] = N[v]$,
\item no F-set if $G$ has open or closed twins,
\end{itemize}
whereas every graph has an L-set.

Accordingly, we call a graph $G$ \emph{TD-admissible} if $G$ has a TD-set and \emph{S-admissible} if $G$ has an S-set for $S \in \{O, I, F\}$.
Similarly, we call a graph $G$ \emph{X-admis\-sible} if $G$ has an X-code for $X \in \codes = \{LD, LTD, OD, OTD, ID, ITD, FD, FTD\}$, and conclude from the previous conditions and the definitions of X-codes which graphs are X-admis\-sible and in which graphs no X-codes exist, see Table \ref{tab_admissible} for a summary. 
\begin{table}[ht]
\begin{center}
\begin{tabular}{ r || c | c || c | c || c | c || c | c }
 X                 & LD   & LTD  & OD   & OTD  & ID & ITD & FD & FTD \\ \hline
 isolated vertices &      &  x   &      &  x   &      &  x   &      &  x   \\
 open twins        &      &      &  x   &  x   &      &      &  x   &  x   \\
 closed twins      &      &      &      &      &  x   &  x   &  x   &  x  \\
\end{tabular}
\end{center}
\caption{
The X-codes that do not exist in graphs having isolated vertices or open or closed twins are marked with an x in the table.}
\label{tab_admissible}
\end{table}



Given an $X \in \codes$ and an X-admis\-sible graph $G$, the X-problem on $G$ is the problem of finding an X-code of minimum cardinality $\x(G)$ in $G$, called its \emph{X-number}. %
Problems of this type have been actively studied in the context of various applications during the last decade, see e.g. the bibliography maintained by Jean and Lobstein \cite{Lobstein_Lib}. 
All X-problems for any $X \in \codes$ have been shown to be NP-hard \cite{CW_ISCO2024,CW_Caldam2025,CHL_2003,CSS_1987,SS_2010}.
Therefore, a typical line to attack such problems is to determine closed formulas for special graphs or to provide bounds on the X-numbers.
We again refer to \cite{Lobstein_Lib} for a literature overview. 

\begin{figure}[h]
\begin{center}
\includegraphics[scale=0.8]{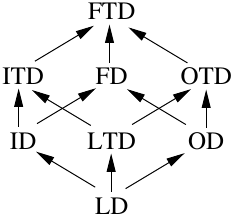}
\caption{The relations between the X-numbers for all $X \in \codes$, where $X' \longrightarrow X$ stands for $\gamma^{X'}(G) \leq \gamma^X(G)$.}
\label{fig_relations}
\end{center}
\end{figure}
A general relation between the X-numbers for all $X \in \codes$ has been established in \cite{CW_Caldam2025},
see Figure \ref{fig_relations}.
We see in particular that $\ld(G)$ is a lower bound for \emph{all} other X-numbers, whereas $\ftd(G)$ is an upper bound for \emph{all} other X-numbers in FTD-admissible graphs. Moreover, it has been shown in \cite{CW_ISCO2024} (resp. \cite{CW_Caldam2025})
that the OD- and OTD-numbers (resp. the FD- and FTD-numbers) of any OTD-admissible (resp. FTD-admissible) graph differ by at most one.

\begin{figure}[!t]
\begin{center}
\includegraphics[scale=1.0]{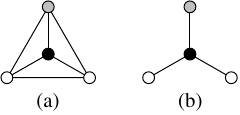}
\caption{Minimum D-sets and TD-sets in two graphs having a universal vertex (the black vertices are universal and form the D-sets, the grey vertices complete the D-sets to the TD-sets), where (a) is the clique $K_4$, (b) is the star $K_{1,3}$.}
\label{fig_universal}
\end{center}
\end{figure}

Moreover, the domination (resp. total-domination) number $\dd(G)$ (resp. $\td(G)$) of a graph $G$, called \emph{D-number} (resp. \emph{TD-number}), is the cardinality of a minimum D-set (resp. TD-set) of $G$.
Both numbers are well-studied in the literature\footnote{In the literature, the domination number of a graph $G$ is typically denoted by $\gamma(G)$ and its total-domination number by $\gamma^t(G)$. However, due to consistency in notation, we prefer to denote them by $\dd(G)$ and $\td(G)$ in this article.}. 
By definition, it is clear that
\begin{equation}\label{eq_D-TD}
  \dd(G) \leq \td(G) \mbox{ holds for any TD-admis\-sible graph $G$}
\end{equation}
and that we have
\begin{equation}\label{ILP_MCF}
  \begin{array}{rl}
   \dd(G) & \!\!\!\!  \leq \x(G) \ \forall X \in \{LD, OD, ID, FD\} \\
   \td(G) & \!\!\!\!  \leq \x(G) \ \forall X \in \{LTD, OTD, ITD, FTD\} \\
\end{array}
\end{equation}
for all X-admissible graphs.
It is well-known that the gap between $\dd(G)$ resp. $\td(G)$ and $\x(G)$ of a graph $G$ can be large.
For instance, if $G$ has a \emph{universal vertex}, that is a vertex $u$ being adjacent to all other vertices of $G$, then $\{u\}$ is a D-set of $G$, whereas $\{u,v\}$ is a TD-set of $G$ for an arbitrary neighbor $v$ of $u$, hence $\dd(G)=1$ and $\td(G)=2$ holds for any graph $G$ having a universal vertex.
Cliques $K_n$ and stars $K_{1,n}$ are clearly graphs with universal vertices, see Figure \ref{fig_universal} for illustration, and provide examples of graphs with a large gap to the X-numbers as results from \cite{ABLW_2022,CW_ISCO2024} show that
$\x(K_n) = n-1 \mbox{ holds for } X \in \{LD, LTD, OD, OTD\}$ 
whereas
$\x(K_{1,n}) = n \mbox{ holds for } X \in \{LD, LTD, ID, ITD\}$
by \cite{ABLW_2022,ABLW_2018-DAM,FL_2023}. 

To the contrary, to the best of our knowledge, the four separation properties have not yet been studied independently, that is not in combination with a domination property to define the eight X-codes with $X \in \codes$.
The aim of this paper is to study these four separation properties.
In Section \ref{sec_complex}, we first show that finding minimum S-sets in a graph is NP-hard for all $S \in \{L, O, I, F\}$.
This motivates to provide different bounds on the cardinality of minimum S-sets, in particular 
in relation to minimum X-codes in X-admissible graphs $G$ (see Section \ref{sec_rel}) and their complements $\bar G$ (see Section \ref{sec_coG}). 
We close with some concluding remarks and lines of future research.

\section{Complexities of computing minimum S-Separating sets}
\label{sec_complex}

Given an S-admissible graph $G$ for some $S \in \{L, O, I, F\}$, the S-problem on $G$ is the problem of finding an S-set of minimum cardinality $\sset(G)$ in $G$, called its \emph{S-number}. In this section, we prove that finding such minimum S-sets, for any $S \in \{L, O, I, F\}$, is NP-complete. The decision version of the problem is stated as follows. 

\defproblem{\minset{S}}{A graph $G$ and a parameter $k$.}{Does there exist an S-set $A$ of $G$ such that $|A| \le k$?}

To prove the NP-completeness of \minset{S}, we reduce the latter from a well-known NP-complete problem called \test (Problem~SP6 in the book by Garey and Johnson~\cite{GJ_1979}) defined as follows.

\defproblem{\test}{A set $U$ of elements called \emph{items}; a collection $\calT$, called a \emph{test collection}, defined as a set of subsets (called \emph{tests}) of $U$ such that for any pair of distinct items there exists a test to which exactly one of the items belongs; and a parameter $\ell$.}{Does there exist a test collection $\calT' \subset \calT$ such that $|\calT'| \le \ell$?}

We next present a general reduction scheme from \test to \minset{S} for $S \in \{L, O, I, F\}$.

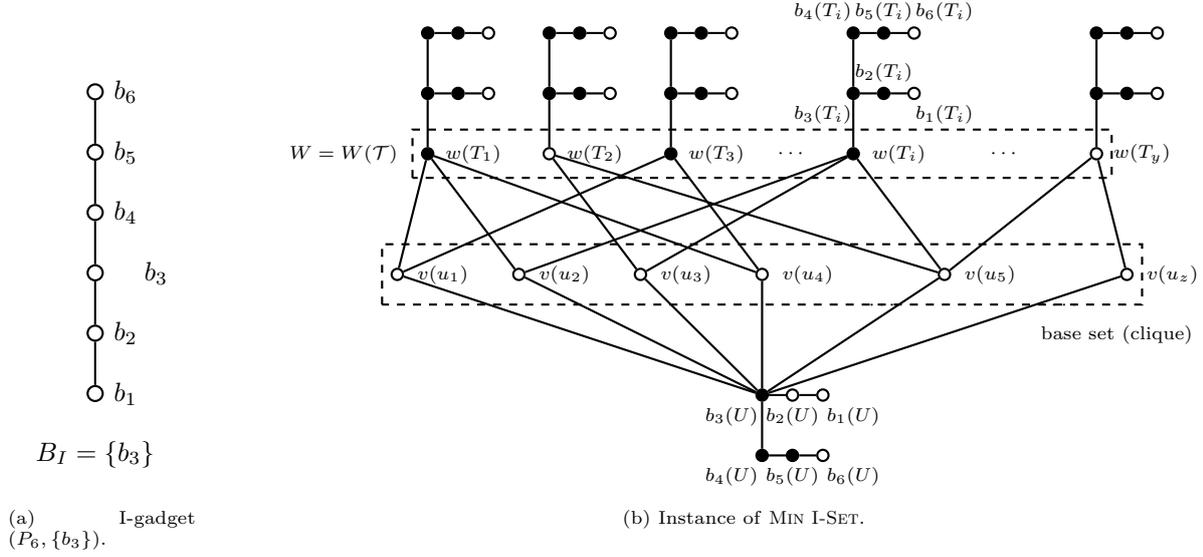
\begin{figure}[!t]
\centering
\begin{subfigure}[t]{0.15\textwidth}
\centering
\begin{tikzpicture}[blacknode/.style={circle, draw=black!, fill=black!, thick, inner sep=2pt},
whitenode/.style={circle, draw=black!, fill=white!, thick, inner sep=2pt},
scale=0.4]


\node[whitenode] (6) at (0,4) {};\node at (1,4) {$b_6$};
\node[whitenode] (5) at (0,2) {};\node at (1,2) {$b_5$};
\node[whitenode] (4) at (0,0) {};\node at (1,0) {$b_4$};
\node[whitenode] (3) at (0,-2) {};\node at (2,-2) {$b_3$};
\node[whitenode] (2) at (0,-4) {};\node at (1,-4) {$b_2$};
\node[whitenode] (1) at (0,-6) {};\node at (1,-6) {$b_1$};

\node () at (0,-8) {$B_I = \{b_3\}$};

\node () at (0,-9) {};

\draw[-, thick, black!] (1) -- (2);
\draw[-, thick, black!] (2) -- (3);
\draw[-, thick, black!] (3) -- (4);
\draw[-, thick, black!] (4) -- (5);
\draw[-, thick, black!] (5) -- (6);


\end{tikzpicture}
\caption{I-gadget $(P_6, \{b_3\})$.} \label{fig:C-gadget}
\end{subfigure}
\hspace{4mm}
\begin{subfigure}[t]{0.8\textwidth}
\centering
\begin{tikzpicture}[blacknode/.style={circle, draw=black!, fill=black!, thick, inner sep=1.5pt},
whitenode/.style={circle, draw=black!, fill=white!, thick, inner sep=1.5pt},
scale=0.4]
\scriptsize

\node[whitenode] (T1_1) at (3,7) {};
\node[blacknode] (T1_2) at (2,7) {};
\node[blacknode] (T1_3) at (1,7) {};
\node[blacknode] (T1_4) at (1,9) {};
\node[blacknode] (T1_5) at (2,9) {};
\node[whitenode] (T1_6) at (3,9) {};

\node[whitenode] (T2_1) at (7,7) {};
\node[blacknode] (T2_2) at (6,7) {};
\node[blacknode] (T2_3) at (5,7) {};
\node[blacknode] (T2_4) at (5,9) {};
\node[blacknode] (T2_5) at (6,9) {};
\node[whitenode] (T2_6) at (7,9) {};

\node[whitenode] (T3_1) at (11,7) {};
\node[blacknode] (T3_2) at (10,7) {};
\node[blacknode] (T3_3) at (9,7) {};
\node[blacknode] (T3_4) at (9,9) {};
\node[blacknode] (T3_5) at (10,9) {};
\node[whitenode] (T3_6) at (11,9) {};

\node[whitenode] (T4_1) at (17,7) {};\node at (18,6.3) {$b_1(T_i)$};
\node[blacknode] (T4_2) at (16,7) {};\node at (16,7.7) {$b_2(T_i)$};
\node[blacknode] (T4_3) at (15,7) {};\node at (14,6.3) {$b_3(T_i)$};
\node[blacknode] (T4_4) at (15,9) {};\node at (14,9.7) {$b_4(T_i)$};
\node[blacknode] (T4_5) at (16,9) {};\node at (16,9.7) {$b_5(T_i)$};
\node[whitenode] (T4_6) at (17,9) {};\node at (18,9.7) {$b_6(T_i)$};

\node[whitenode] (T5_1) at (25,7) {};
\node[blacknode] (T5_2) at (24,7) {};
\node[blacknode] (T5_3) at (23,7) {};
\node[blacknode] (T5_4) at (23,9) {};
\node[blacknode] (T5_5) at (24,9) {};
\node[whitenode] (T5_6) at (25,9) {};

\node[blacknode] (1') at (1,5) {};\node at (2.5,5) {$w(T_1)$};0
\node[whitenode] (2') at (5,5) {};\node at (6.5,5) {$w(T_2)$};
\node[blacknode] (3') at (9,5) {};\node at (10.5,5) {$w(T_3)$};
\node at (13,5) {$\cdots$};
\node[blacknode] (4') at (15,5) {};\node at (16.5,5) {$w(T_i)$};
\node at (20,5) {$\cdots$};
\node[whitenode] (5') at (23,5) {};\node at (24.5,5) {$w(T_y)$};

\node[whitenode] (1) at (0,1) {};\node at (1.5,1) {$v(u_1)$};
\node[whitenode] (2) at (4,1) {};\node at (5.5,1) {$v(u_2)$};
\node[whitenode] (3) at (8,1) {};\node at (9.5,1) {$v(u_3)$};
\node[whitenode] (4) at (12,1) {};\node at (13.5,1) {$v(u_4)$};
\node at (16,0) {$\cdots$};
\node[whitenode] (5) at (18,1) {};\node at (19.5,1) {$v(u_5)$};
\node at (22,0) {$\cdots$};
\node[whitenode] (6) at (24,1) {};\node at (25.5,1) {$v(u_z)$};


\node[whitenode] (U_1) at (14,-3) {};\node at (15,-3.7) {$b_1(U)$};
\node[whitenode] (U_2) at (13,-3) {};\node at (13,-3.7) {$b_2(U)$};
\node[blacknode] (U_3) at (12,-3) {};\node at (11,-3.7) {$b_3(U)$};
\node[blacknode] (U_4) at (12,-5) {};\node at (11,-5.7) {$b_4(U)$};
\node[blacknode] (U_5) at (13,-5) {};\node at (13,-5.7) {$b_5(U)$};
\node[whitenode] (U_6) at (14,-5) {};\node at (15,-5.7) {$b_6(U)$};

\draw[dashed, thick] (-0.5,0) rectangle (24.5,2) node[below right = 0.5cm and 3.5cm of 4] {base set (clique)};

\draw[dashed, thick] (0.5,4.2) rectangle (23.5,5.8) node[left = 0.2cm of 1'] {$W = W(\calT)$};

\draw[-, thick, black!] (T1_1) -- (T1_2);
\draw[-, thick, black!] (T1_2) -- (T1_3);
\draw[-, thick, black!] (T1_3) -- (T1_4);
\draw[-, thick, black!] (T1_4) -- (T1_5);
\draw[-, thick, black!] (T1_5) -- (T1_6);

\draw[-, thick, black!] (T2_1) -- (T2_2);
\draw[-, thick, black!] (T2_2) -- (T2_3);
\draw[-, thick, black!] (T2_3) -- (T2_4);
\draw[-, thick, black!] (T2_4) -- (T2_5);
\draw[-, thick, black!] (T2_5) -- (T2_6);

\draw[-, thick, black!] (T3_1) -- (T3_2);
\draw[-, thick, black!] (T3_2) -- (T3_3);
\draw[-, thick, black!] (T3_3) -- (T3_4);
\draw[-, thick, black!] (T3_4) -- (T3_5);
\draw[-, thick, black!] (T3_5) -- (T3_6);

\draw[-, thick, black!] (T4_1) -- (T4_2);
\draw[-, thick, black!] (T4_2) -- (T4_3);
\draw[-, thick, black!] (T4_3) -- (T4_4);
\draw[-, thick, black!] (T4_4) -- (T4_5);
\draw[-, thick, black!] (T4_5) -- (T4_6);

\draw[-, thick, black!] (T5_1) -- (T5_2);
\draw[-, thick, black!] (T5_2) -- (T5_3);
\draw[-, thick, black!] (T5_3) -- (T5_4);
\draw[-, thick, black!] (T5_4) -- (T5_5);
\draw[-, thick, black!] (T5_5) -- (T5_6);

\draw[-, thick, black!] (T1_3) -- (1');
\draw[-, thick, black!] (T2_3) -- (2');
\draw[-, thick, black!] (T3_3) -- (3');
\draw[-, thick, black!] (T4_3) -- (4');
\draw[-, thick, black!] (T5_3) -- (5');

\draw[-, thick, black!] (U_1) -- (U_2);
\draw[-, thick, black!] (U_2) -- (U_3);
\draw[-, thick, black!] (U_3) -- (U_4);
\draw[-, thick, black!] (U_4) -- (U_5);
\draw[-, thick, black!] (U_5) -- (U_6);

\draw[-, thick, black!] (1) -- (U_3);
\draw[-, thick, black!] (2) -- (U_3);
\draw[-, thick, black!] (3) -- (U_3);
\draw[-, thick, black!] (4) -- (U_3);
\draw[-, thick, black!] (5) -- (U_3);
\draw[-, thick, black!] (6) -- (U_3);

\draw[-, thick, black!] (1) -- (1');
\draw[-, thick, black!] (1) -- (3');
\draw[-, thick, black!] (2) -- (1');
\draw[-, thick, black!] (2) -- (4');
\draw[-, thick, black!] (3) -- (2');
\draw[-, thick, black!] (3) -- (4');
\draw[-, thick, black!] (4) -- (1');
\draw[-, thick, black!] (4) -- (3');
\draw[-, thick, black!] (5) -- (2');
\draw[-, thick, black!] (5) -- (4');
\draw[-, thick, black!] (5) -- (5');
\draw[-, thick, black!] (6) -- (5');

\end{tikzpicture}
\caption{Instance of \minset{I}.} \label{fig:NP_C_all}
\end{subfigure}
\caption[]{Illustration of the graph $G_I$ constructed by Reduction~\ref{red:general}. Here $U = \{u_1, u_2, \ldots u_z\}$ and $\calT = \{T_1, T_2, \ldots , T_y\}$. The set of black vertices in Figure~\ref{fig:NP_C_all} is an example of an I-set of $G_I$.}
\label{fig:NP_C}
\end{figure}

\begin{reduction}[General scheme] \label{red:general}
The reduction takes as input a set $U$ of items, a test collection $\calT$ of $U$ and a parameter $\ell$. It then creates a graph $G_S$ by the following steps (see Figure~\ref{fig:NP_C}).

\begin{enumerate}[leftmargin=16pt, itemsep=0pt]
\item[(1)] The reduction sets a positive integer $r = r(S)$ such that for each $u \in U$, it introduces a set $V(S,u) = \{v^1(u), v^2(u) \ldots , v^r(u)\}$ of vertices of $G_S$. Let $M = M(S) = \cup_{u \in U} V(S,u)$. It also introduces a set $R = R(S)$ of vertices of size $\calO (\ell)$. Then the set $Q = Q(S) = M \cup R$ is called a \emph{base set} of $G_S$. Note that $|Q| = r|U| + |R|$. Depending on the separating property $S$, the reduction then introduces edges between pairs of vertices in $Q$.
\item[(2)] For each $T \in \calT$, the reduction introduces a vertex $w = w(\calT)$ of $G_S$. It then introduces edges such that $v^i(u) w(T) \in E(G)$ for all $i \in \{1,2, \ldots , r\}$ if and only if $u \in T$. Let $W = W(T) = \{w(T) : T \in \calT\}$.
\item[(3)] We define an \emph{S-gadget} to be a tuple $(H_S,B_S)$ such that $H_S$ is a graph and $B_S$ is a set of some \emph{designated vertices} of $H_S$. For each test $T \in \calT$, the reduction introduces an S-gadget $(H_S(T),B_S(T))$ and adds an edge between $w(T)$ and each $b \in B_S(T)$. In addition, it introduces another S-gadget $(H_S(U), b_S(U))$ and edges between each $b \in B_S(U)$ and $v^i(u) \in M$ for each $i \in \{1, 2, \ldots , \ell+1\}$ and $u \in U$. See Figure~\ref{fig:C-gadget} for an example when $S = I$.
\item[(4)] Finally, the reduction sets $k = \ell + p |\calT| + q$, where $p = p(S)$ and $q = q(S)$ are integers.
\end{enumerate}
\end{reduction}
Notice that Reduction~\ref{red:general} can be executed in polynomial-time in the input size of \test.

We next illustrate how we adapt Reduction~\ref{red:general} to treat closed separation, that is for $S = I$. In this case, in Step (1) of Reduction~\ref{red:general}, the reduction sets $r=r(I) =1$ and so $V(I,u) = \{v(u)\}$ for each $u \in U$. Moreover, it sets $R = R(I) = \emptyset$. The base set $Q = Q(I) = M(I) = \cup_{u \in U} V(I,u)$ is therefore of order $|U|$ and the reduction introduces edges between all pairs of vertices of $Q$ making $G_I[Q]$ a complete graph. The I-gadget in Step (3) of Reduction~\ref{red:general} is $(H_I,B_I)$, where $H_I$ is isomorphic to $P_6$ and $B_I = \{b_3\}$ as shown in Figure~\ref{fig:C-gadget}. Finally, the reduction sets $p = p(I) = 4$ and $q = q(I) = 3$ in Step (4) of Reduction~\ref{red:general}, that is, it sets $k = \ell + 4 |\calT| + 3$.

We use the same vertex names as mentioned in Figure~\ref{fig:NP_C}. Moreover, for notational convenience, let $V_U = V(H_I(U))$ and let $V_T = V(H_I(T))$ for each $T \in \calT$. The following lemma is helpful in proving the NP-hardness of \minset{I} next.


\begin{lemma} \label{lem:NP_C}
Let $A$ be an I-set of $G_I$. Then we have $\left |A \cap \Big ( V_U \cup \bigcup_{T \in \calT} V_T \Big ) \right | \ge 4|\calT|+3$.
\end{lemma}

With the help of Lemma \ref{lem:NP_C}, we can prove the following.

\begin{theorem}\label{thm_I-hard}
\minset{I} is NP-complete.
\end{theorem}

\begin{proof}[Proof (sketch)]
The problem clearly belongs to the class NP. Hence, to prove the result, we show that \minset{I} is also NP-hard. We do so by showing that an input $(U,\calT)$ to \test is a YES-instance of \test if and only if the input $(G_I,k)$ constructed in Reduction~\ref{red:general} is a YES-instance of \minset{I}. For the necessary part, given a test collection $\calT'$ of $U$ such that $|\calT'| \le \ell$, the idea is to build an I-set of $G_I$ by including the corresponding vertices $w(T)$ for $T \in \calT'$ and those marked in black in Figure~\ref{fig:NP_C_all} in the $(|\calT|+1)$ induced subgraphs $H_I$ of $G_I$. On the other hand, for the suffiency part, if there exists an I-set of $G_I$ of cardinality at most $\ell + 4 |\calT| + 3$, then Lemma~\ref{lem:NP_C} ensures that pairs of vertices in $Q$ are closed-separated by at most $\ell$ vertices in $W = W(\calT)$. Hence, the set $\calT' \subset \calT$ of their corresponding tests is a test collection of $U$ with $|\calT'| \le \ell$. This proves the result.
\end{proof}

As in the cases for $S =I$, Reduction~\ref{red:general} can also be adapted for the cases when $S \in \{L,O,F\}$ to prove the next theorem. Due to a limitation to the number of pages, we leave the rest of the proofs in the appendix.

\begin{theorem} \label{thm:NP_main}
\minset{S} is NP-complete for all $S \in \{L,I,O,F\}$.
\end{theorem}


\section{From S-sets to X-codes}\label{sec_rel}

The hardness results from the previous section motivate to study the relations of S-numbers and some bounds for their values. Given the general relations between the X-numbers for all $X \in \codes$ in Figure \ref{fig_relations}, we immediately conclude that
$
\lset(G) \leq \oset(G),\iset(G) \leq \fset(G)
$ 
holds for all accordingly \sadmis ~graphs $G$ as well as 
\begin{equation}\label{eq_Set-leq-Codes}
   \sset(G) \leq \sd(G) \leq \std(G) \ \forall S \in \{L, O, I, F\} \\
\end{equation}
for all accordingly X-admissible graphs. Next,
we establish upper bounds on the X-number of a graph in terms of its S-number. 
We start with the relation of the S-number of a graph and its SD-number: 

\begin{theorem} \label{thm_S-2-SD-code}
  For any $S \in \{L, O, I, F\}$ and an \sadmis ~graph $G$, we have $\sd(G) \leq \sset(G)+1$.
\end{theorem}

We next wonder whether the same relation holds for the S-number of a graph and its STD-number.
On the one hand, we can show:

\begin{theorem} \label{thm_S-2-STD-code1}
  For $S \in \{O, F\}$ and an \sadmis ~graph $G$, we have $\std(G) \leq \sset(G)+1$.
\end{theorem}

On the other hand, we have:

\begin{theorem} \label{thm_S-2-STD-code2}
  For $S \in \{L, I\}$ and an \sadmis ~graph $G$, we have $\std(G) \leq 2 \sset(G)$.
\end{theorem}

\begin{figure}[!t]
\begin{center}
\includegraphics[scale=1.0]{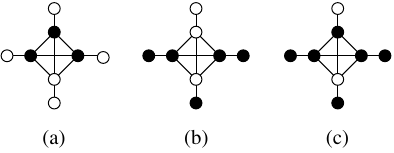}
\caption{Minimum S-sets in the thin headless spider $H_4$ (the black vertices belong to the S-sets), where (a) is both an L-set and an O-set, (b) is an I-set, (c) is an F-set.}
\label{fig_thin-spiders}
\end{center}
\end{figure}

For illustration, we are going to determine the S-numbers of a family of graphs without open or closed twins for which all X-numbers are known from the literature. 
A \emph{thin headless spider} is a graph $H_k=(Q \cup S, E)$ where $Q = \{q_1, \ldots, q_k\}$ induces a clique, $S = \{s_1, \ldots, s_k\}$ a stable set, and $q_i$ is adjacent to $s_j$ if and only if $i = j$.
We have, for example, $H_2 = P_4$ and the thin headless spider $H_4$ is depicted in Figure~\ref{fig_thin-spiders}.

We can show the following for the S-numbers of thin headless spiders:

\begin{theorem} \label{thm_thin-spiders}
For a thin headless spider $H_k=(Q \cup S, E)$ with $k \geq 4$, we have 
\[
\sset(H_k) = \left\{
\begin{array}{rl}
  k-1   & \mbox{if } S \in \{L,O\}  \\
  k+1   & \mbox{if } S = I,  \\
  2k-2  & \mbox{if } S = F.  \\
\end{array}
\right.
\]
\end{theorem}

Comparing those S-numbers with the X-numbers of thin headless spiders known from the literature 
as summarized in Table \ref{tab_thin-spiders}, we clearly see that $\sset(H_k)$ and $\sd(H_k)$ differ indeed by at most 1 for all $S \in \{L, O, I, F\}$, that also $\sset(H_k)$ and $\std(H_k)$ differ by only 1 for $S \in \{O, F\}$, whereas the gap between $\iset(H_k)$ and $\itd(H_k)$ is indeed large.

\begin{table}[h]
\begin{center}
\begin{tabular}{ r || c | c | c | c  }
 S            & L                & O                  & I                 &  F \\ \hline
 $\sset(H_k)$ & $k-1$            & $k-1$              & $k+1$             &  $2k-2$    \\
 $\sd(H_k)$   & $k$ by \cite{ABLW_2022} & $k$ by \cite{CW_ISCO2024} & $k+1$ by \cite{ABW_2016} &  $2k-2$ by \cite{CW_2024-FS}   \\
 $\std(H_k)$  & $k$ by \cite{ABLW_2022} & $k$ by \cite{ABLW_2022}   & $2k-1$ by \cite{FL_2023} &  $2k-1$ by \cite{CW_2024-FS}  \\
\end{tabular}
\end{center}
\caption{
The S- and X-numbers of thin headless spiders $H_k$ for $k \geq 4$ and all $S \in \{L, O, I, F\}$.}
\label{tab_thin-spiders}
\end{table}

\section{S-sets, X-codes and complementation}\label{sec_coG}

Given an S-admissible graph $G$ for some $S \in \{L, O, I, F\}$, we first address the question concerning the S-numbers of $G$ and its complement $\bar G$. We can show:\\

\begin{figure}[!t]
\begin{center}
\includegraphics[scale=0.9]{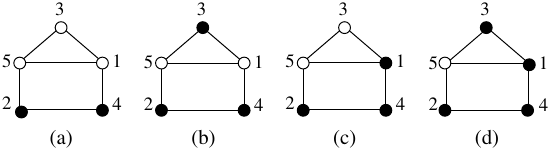}
\caption{Minimum S-sets in a graph (the black vertices belong to the S-sets), where (a) is an L-set, (b) is an O-set (and also a D-set, hence also an OD-code), (c) is an I-set (and also a TD-set, hence also both an ID- and an ITD-code), (d) is an F-set (and also a TD-set, hence both an FD- and an FTD-code).}
\label{fig_exp_coS-sets}
\end{center}
\end{figure}

\begin{theorem} \label{thm_coS-numbers}
  For any graph $G$
  \begin{itemize}
    \itemsep -3pt
  \item[(a)] we have $\lset(G) = \lset(\bar G)$;
  \item[(b)] without closed twins, $\iset(G) = \oset(\bar G)$ holds;
  \item[(c)] without open twins, we have $\oset(G) = \iset(\bar G)$;
  \item[(d)] without twins, $\fset(G) = \fset(\bar G)$ follows.
  \end{itemize}
\end{theorem}
For illustration, consider the minimum S-sets of the 5-path shown in Figure \ref{fig_exp_X-codes} and compare them with the minimum S-sets of the complement of the 5-path, the house, shown in Figure \ref{fig_exp_coS-sets}.

Combining the above results with the findings from Section \ref{sec_rel}, we further conclude:

\begin{corollary}\label{cor_coSD-numbers}
  For any graph $G$
  \begin{itemize}
    \itemsep -3pt
  \item[(a)] $\ld(G)$ and $\ld(\bar G)$ differ by at most one;
  \item[(b)] without closed twins, $\id(G)$ and $\od(\bar G)$ differ by at most one;
  \item[(c)] without open twins, $\od(G)$ and $\id(\bar G)$ differ by at most one;
  \item[(d)] without twins, $\fd(G)$ and $\fd(\bar G)$ as well as $\ftd(G)$ and $\ftd(\bar G)$ differ by at most one.
  \end{itemize}
\end{corollary}

For illustration, we consider the complements of thin headless spiders $H_k$, called \emph{thick headless spiders}, which are graphs $\bar{H}_k=(Q \cup S, E)$ where $Q = \{q_1, \ldots, q_k\}$ induces a stable set, $S = \{s_1, \ldots, s_k\}$ a clique, and $q_i$ is adjacent to $s_j$ if and only if $i \neq j$.

Also for thick headless spiders, all X-numbers are known from the literature and we obtain their S-numbers by combining the findings of Theorem \ref{thm_thin-spiders} and Theorem \ref{thm_coS-numbers}. 
Comparing those S-numbers with the X-numbers as summarized in Table \ref{tab_thick-spiders}, we see that
$\sset(\bar{H}_k)$ equals $\sd(\bar{H}_k)$ and $\std(\bar{H}_k)$ for all $S \in \{L, O, F\}$, that
$\iset(\bar{H}_k)$ and $\id(\bar{H}_k)$ differ by only 1 but 
$\iset(\bar{H}_k)$ and $\itd(\bar{H}_k)$ by 2.

\begin{table}[ht]
\begin{center}
\begin{tabular}{ r || c | c | c | c  }
 S                  & L                  & I                     & O                    &  F \\ \hline
 $\sset(\bar{H}_k)$ & $k-1$              & $k-1$                 & $k+1$                &  $2k-2$    \\
 $\sd(\bar{H}_k)$   & $k-1$ by \cite{ABLW_2022} & $k$ by \cite{ABW_2016}       & $k+1$ by \cite{CW_ISCO2024} &  $2k-2$ by \cite{CW_2024-FS}   \\
 $\std(\bar{H}_k)$  & $k-1$ by \cite{ABLW_2022} & $k+1$ by \cite{CW_2024-FS} & $k+1$ by \cite{ABLW_2022}  &  $2k-2$ by \cite{CW_2024-FS}  \\
\end{tabular}
\end{center}
\caption{
The S- and X-numbers of thick headless spiders $\bar{H}_k$ for $k \geq 4$ and all $S \in \{L, O, I, F\}$.}
\label{tab_thick-spiders}
\end{table}

\section{Conclusion}

In this paper, we studied the four separation properties location, closed-separation, open-separation and full-separation. 
We showed that finding minimum separating sets in a graph is NP-hard, that is, as hard as the related problems for (total-)dominating sets or X-codes.
Then we focused on the interplay of the studied separation and domination properties, revealing that separation is almost domination due to 
$$
\sset(G) \leq \sd(G) \leq \sset(G) + 1 \ \forall S \in \{L, O, I, F\}
$$
and all SD-admissible graphs $G$. Moreover, we showed
$$
\sset(G) \leq \sd(G) \leq \std(G) \leq \sset(G) + 1 \mbox{ for } S \in \{O, F\}
$$
which reproves also the according relations of SD- and STD-numbers for $S \in \{O, F\}$ established in \cite{CW_ISCO2024,CW_Caldam2025}
whereas 
$$
\sset(G) \leq \std(G) \leq 2 \sset(G)  \mbox{ for } S \in \{L, I\}
$$
shows that the two S-numbers can substantially differ from the according STD-numbers.
The latter result for $S=I$ reproves a relation between ID- and ITD-numbers established in \cite{FL_2023}. 

Moreover, studying the interplay of separation and complementation revealed that location and full-separation are the same on a graph and its complement, whereas closed-separation in a graph corresponds to open-separation in its complement, with according implications for the X-numbers of a graph and its complement.

Our lines of future research include to provide tight examples for all established bounds and relations as well as to study the complexity of the problems to decide whether S- and X-numbers that differ by at most one are in fact equal for a given graph.
Seen that deciding whether $\sd(G)$ equals $\std(G)$ for $S \in \{O, F\}$ was proven to be NP-hard in \cite{CW_ISCO2024,CW_Caldam2025}
we conjecture that, given an SD-admissible graph $G$, it is also NP-hard to decide whether $\sset(G)$ equals $\sd(G)$ for all $S \in \{L, O, I,  F\}$.





\bibliographystyle{elsarticle-num}
\bibliography{idproblembib}


\newpage


\begin{center}
\huge
\textbf{Appendix}
\end{center}

\vspace{1cm}

\noindent In the appendix, we present the proofs of the results that have been omitted due to restrictions to the number of pages.

\section*{Proofs omitted in Section~\ref{sec_complex}}

We use the same vertex names as mentioned in Figure~\ref{fig:NP_C}. For an I-gadget (whose vertices are named as $b_j$ in Figure~\ref{fig:C-gadget}), the vertices of the graph $H_I(T)$ for some $T \in \calT$ or the graph $H_I(U)$ are denoted as $b_j(T)$ or $b_j(U)$, respectively. Moreover, for notational convenience, let $V_U = V(H_I(U))$ and let $V_T = V(H_I(T))$ for each $T \in \calT$.



\begin{proof}[Proof of Lemma~\ref{lem:NP_C}]
Let $A$ be an I-set of $G_I$. For each $T \in \calT$, we must have $b_3(T), b_4(T) \in A$ in order for the pairs $b_1(T), b_2(T)$ and $b_5(T), b_6(T)$, respectively, to be I-separated. By the same argument, we also have $b_3(U), b_4(U) \in A$. Moreover, among all the $(|\calT|+1)$ number of I-gadgets $(H_I,B_I)$ introduced by Reduction~\ref{red:general}, there can exist at most one support vertex and its adjacent leaf of one such $H_I$ which do not belong to $A$ (or else, there would be at least two leaves with empty closed neighborhoods in $A$ contradicting the fact that $A$ is an I-set of $G_I$). Thus, $A$ contains at least $4$ vertices from all $(|\calT|+1)$ number of $H_I$ introduced by Reduction~\ref{red:general} except possibly for one $H_I$ from which $A$ contains $3$ vertices. This proves the result.
\end{proof}

\begin{proof}[Proof of Theorem~\ref{thm_I-hard}]
The problem clearly belongs to the class NP. Hence, to prove the result, we show that \minset{I} is also NP-hard. We do so by showing that an input $(U,\calT)$ to \test is a YES-instance of \test if and only if the input $(G_I,k)$ constructed in Reduction~\ref{red:general} is a YES-instance of \minset{I}. 

($\implies$): Let $(U,\calT)$ be a YES-instance of \test. In other words, there exists a test collection $\calT'$ of $U$ such that $|\calT'| \le \ell$. We then construct a set $A \subset V(G_I)$ as follows. For each $T \in \calT'$, let $w(T) \in A$. Also, for all $T \in \calT$, let $b_2(T), b_3(T), b_4(T), b_5(T) \in A$ and let $b_3(U), b_4(U), b_5(U) \in A$. This implies that $|A| = \ell + 4 |\calT|+3$. Since $\calT'$ is a test collection of $U$, by the construction of $G_I$, the set $A \cap W$ I-separates each pair of vertices in $Q$. 
Moreover, it can be verified that each vertex of $G_I$ has a unique closed neighborhood in A (with only $b_1(U)$ having an empty closed neighborhood in $A$). This makes $A$ a I-set of $G_I$ and hence, proves the necessary part of the result.

($\impliedby$): Let $(G_I,k)$ be a YES-instance of \minset{I}. In other words, there exists an I-set $A$ of $G_I$ such that $|A| \le k = \ell + 4 |\calT| + 3$. Now, using Lemma~\ref{lem:NP_C}, we must have $|A \cap W| \le \ell$. We now create a subcollection $\calT'$ of $\calT$ by by letting $T \in \calT'$ if $w(T) \in A$. Then, we have $|\calT'| \le \ell$. Since $G_S[Q \cup \{b_4(U)\}]$ is a complete graph, the neighborhoods of the vertices $Q$ can differ only by their neighborhoods in $W$. Thus, pairs of vertices in $Q$ are I-separated by the set $A \cap W$. Hence, by the construction of $G_I$, 
the collection $\calT'$ is a test collection of $U$. This proves the sufficiency part of the result and hence, proves the theorem.
\end{proof}

Toward the proof of Theorem~\ref{thm:NP_main}, in what follows, we adapt Reduction~\ref{red:general} to the cases where the separating property $S \in \{O,F,L\}$. For the case $S = O$, however, the result holds more easily by the following corollary. Since Lemma \ref{lemA_S_coG} shows that for any graph $G$, we have $\hyp_I(G) = \hyp_O(\overline G)$, it follows that any I-set of a graph $G$ is an O-set of $\overline G$ and vice-versa. Hence determining a minimum O-set of $G$ is equivalent to determining a minimum I-set of $\overline G$. The NP-completeness of \minset{I} thus immediately implies:
\begin{corollary} \label{cor:NP_O}
\minset{O} is NP-complete.
\end{corollary}

\subsection*{Full-separating sets}
We next adapt Reduction~\ref{red:general} for the case $S = F$. In this case, in Step (1) of Reduction~\ref{red:general}, the reduction sets $r=r(F) =1$ and so $V(F,u) = \{v(u)\}$ for each $u \in U$. It also sets $R = R(F) = \emptyset$. The base set $Q = Q(F) = M(F) = \cup_{u \in U} V(F,u)$ is therefore of order $|U|$ and the reduction introduces edges between all pairs of vertices of $Q$ making $G_S[Q]$ a complete graph. The F-gadget in Step (3) of Reduction~\ref{red:general} is $(H_F,B_F)$, where $H_F$ is isomorphic to the graph in Figure~\ref{fig:C-gadget} and  $B_F= \{b_5\}$. Finally, the reduction sets $p = p(F) = 12$ and $q = q(F) = 11$ in Step (4) of Reduction~\ref{red:general}, that is, it sets $k = \ell + 12 |\calT| + 11$.

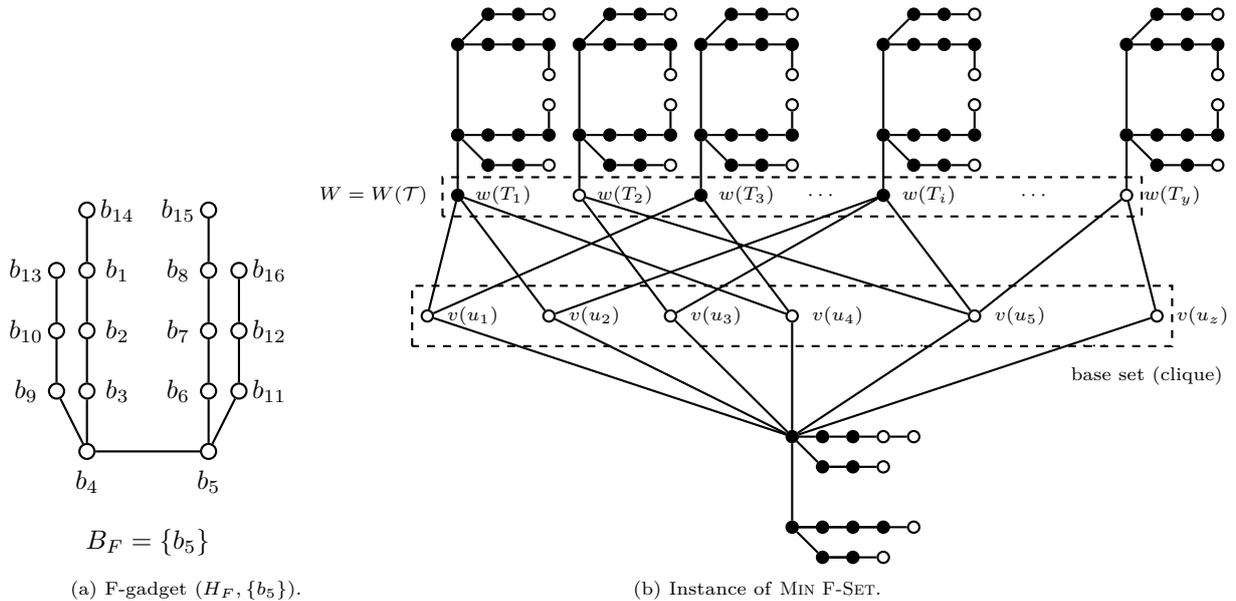
\begin{figure}[!t]
\centering
\begin{subfigure}[t]{0.2\textwidth}
\centering
\begin{tikzpicture}[blacknode/.style={circle, draw=black!, fill=black!, thick, inner sep=2pt},
whitenode/.style={circle, draw=black!, fill=white!, thick, inner sep=2pt},
scale=0.4]
\hspace{-0.8cm}


\node[whitenode] (1) at (0,8) {};\node at (1,8) {$b_{14}$};
\node[whitenode] (2) at (0,6) {};\node at (1,6) {$b_1$};
\node[whitenode] (3) at (0,4) {};\node at (1,4) {$b_2$};
\node[whitenode] (4) at (0,2) {};\node at (1,2) {$b_3$};
\node[whitenode] (5) at (0,0) {};\node at (0,-1) {$b_4$};
\node[whitenode] (6) at (4,0) {};\node at (4,-1) {$b_5$};
\node[whitenode] (7) at (4,2) {};\node at (3,2) {$b_6$};
\node[whitenode] (8) at (4,4) {};\node at (3,4) {$b_7$};
\node[whitenode] (9) at (4,6) {};\node at (3,6) {$b_8$};
\node[whitenode] (10) at (4,8) {};\node at (3,8) {$b_{15}$};
\node[whitenode] (11) at (-1,2) {};\node at (-2,2) {$b_9$};
\node[whitenode] (12) at (-1,4) {};\node at (-2,4) {$b_{10}$};
\node[whitenode] (13) at (-1,6) {};\node at (-2,6) {$b_{13}$};
\node[whitenode] (14) at (5,2) {};\node at (6,2) {$b_{11}$};
\node[whitenode] (15) at (5,4) {};\node at (6,4) {$b_{12}$};
\node[whitenode] (16) at (5,6) {};\node at (6,6) {$b_{16}$};

\node () at (2,-3) {$B_F = \{b_5\}$};


\draw[-, thick, black!] (1) -- (2);
\draw[-, thick, black!] (2) -- (3);
\draw[-, thick, black!] (3) -- (4);
\draw[-, thick, black!] (4) -- (5);
\draw[-, thick, black!] (5) -- (6);
\draw[-, thick, black!] (6) -- (7);
\draw[-, thick, black!] (7) -- (8);
\draw[-, thick, black!] (8) -- (9);
\draw[-, thick, black!] (9) -- (10);
\draw[-, thick, black!] (5) -- (11);
\draw[-, thick, black!] (11) -- (12);
\draw[-, thick, black!] (12) -- (13);
\draw[-, thick, black!] (6) -- (14);
\draw[-, thick, black!] (14) -- (15);
\draw[-, thick, black!] (15) -- (16);


\end{tikzpicture}
\caption{F-gadget $(H_F, \{b_5\})$.} \label{fig:F-gadget}
\end{subfigure}
\begin{subfigure}[t]{0.7\textwidth}
\centering
\begin{tikzpicture}[blacknode/.style={circle, draw=black!, fill=black!, thick, inner sep=1.5pt},
whitenode/.style={circle, draw=black!, fill=white!, thick, inner sep=1.5pt},
scale=0.4]
\scriptsize

\hspace{-0.1cm}

\node[blacknode] (T1_1) at (3,7) {};
\node[blacknode] (T1_2) at (2,7) {};
\node[blacknode] (T1_3) at (1,7) {};
\node[blacknode] (T1_4) at (1,10) {};
\node[blacknode] (T1_5) at (2,10) {};
\node[blacknode] (T1_6) at (3,10) {};
\node[blacknode] (T1_0) at (4,7) {};
\node[blacknode] (T1_7) at (4,10) {};

\node[blacknode] (T1_8) at (2,11) {};
\node[blacknode] (T1_9) at (3,11) {};
\node[blacknode] (T1_10) at (2,6) {};
\node[blacknode] (T1_11) at (3,6) {};
\node[whitenode] (T1_12) at (4,11) {};
\node[whitenode] (T1_13) at (4,6) {};
\node[whitenode] (T1_14) at (4,9) {};
\node[whitenode] (T1_15) at (4,8) {};

\node[blacknode] (T2_1) at (7,7) {};
\node[blacknode] (T2_2) at (6,7) {};
\node[blacknode] (T2_3) at (5,7) {};
\node[blacknode] (T2_4) at (5,10) {};
\node[blacknode] (T2_5) at (6,10) {};
\node[blacknode] (T2_6) at (7,10) {};
\node[blacknode] (T2_0) at (8,7) {};
\node[blacknode] (T2_7) at (8,10) {};

\node[blacknode] (T2_8) at (6,11) {};
\node[blacknode] (T2_9) at (7,11) {};
\node[blacknode] (T2_10) at (6,6) {};
\node[blacknode] (T2_11) at (7,6) {};
\node[whitenode] (T2_12) at (8,11) {};
\node[whitenode] (T2_13) at (8,6) {};
\node[whitenode] (T2_14) at (8,9) {};
\node[whitenode] (T2_15) at (8,8) {};

\node[blacknode] (T3_1) at (11,7) {};
\node[blacknode] (T3_2) at (10,7) {};
\node[blacknode] (T3_3) at (9,7) {};
\node[blacknode] (T3_4) at (9,10) {};
\node[blacknode] (T3_5) at (10,10) {};
\node[blacknode] (T3_6) at (11,10) {};
\node[blacknode] (T3_0) at (12,7) {};
\node[blacknode] (T3_7) at (12,10) {};

\node[blacknode] (T3_8) at (10,11) {};
\node[blacknode] (T3_9) at (11,11) {};
\node[blacknode] (T3_10) at (10,6) {};
\node[blacknode] (T3_11) at (11,6) {};
\node[whitenode] (T3_12) at (12,11) {};
\node[whitenode] (T3_13) at (12,6) {};
\node[whitenode] (T3_14) at (12,9) {};
\node[whitenode] (T3_15) at (12,8) {};

\node[blacknode] (T4_1) at (17,7) {};
\node[blacknode] (T4_2) at (16,7) {};
\node[blacknode] (T4_3) at (15,7) {};
\node[blacknode] (T4_4) at (15,10) {};
\node[blacknode] (T4_5) at (16,10) {};
\node[blacknode] (T4_6) at (17,10) {};
\node[blacknode] (T4_0) at (18,7) {};
\node[blacknode] (T4_7) at (18,10) {};

\node[blacknode] (T4_8) at (16,11) {};
\node[blacknode] (T4_9) at (17,11) {};
\node[blacknode] (T4_10) at (16,6) {};
\node[blacknode] (T4_11) at (17,6) {};
\node[whitenode] (T4_12) at (18,11) {};
\node[whitenode] (T4_13) at (18,6) {};
\node[whitenode] (T4_14) at (18,9) {};
\node[whitenode] (T4_15) at (18,8) {};

\node[blacknode] (T5_1) at (25,7) {};
\node[blacknode] (T5_2) at (24,7) {};
\node[blacknode] (T5_3) at (23,7) {};
\node[blacknode] (T5_4) at (23,10) {};
\node[blacknode] (T5_5) at (24,10) {};
\node[blacknode] (T5_6) at (25,10) {};
\node[blacknode] (T5_0) at (26,7) {};
\node[blacknode] (T5_7) at (26,10) {};

\node[blacknode] (T5_8) at (24,11) {};
\node[blacknode] (T5_9) at (25,11) {};
\node[blacknode] (T5_10) at (24,6) {};
\node[blacknode] (T5_11) at (25,6) {};
\node[whitenode] (T5_12) at (26,11) {};
\node[whitenode] (T5_13) at (26,6) {};
\node[whitenode] (T5_14) at (26,9) {};
\node[whitenode] (T5_15) at (26,8) {};

\node[blacknode] (1') at (1,5) {};\node at (2.5,5) {$w(T_1)$};0
\node[whitenode] (2') at (5,5) {};\node at (6.5,5) {$w(T_2)$};
\node[blacknode] (3') at (9,5) {};\node at (10.5,5) {$w(T_3)$};
\node at (13,5) {$\cdots$};
\node[blacknode] (4') at (15,5) {};\node at (16.5,5) {$w(T_i)$};
\node at (20,5) {$\cdots$};
\node[whitenode] (5') at (23,5) {};\node at (24.5,5) {$w(T_y)$};

\node[whitenode] (1) at (0,1) {};\node at (1.5,1) {$v(u_1)$};
\node[whitenode] (2) at (4,1) {};\node at (5.5,1) {$v(u_2)$};
\node[whitenode] (3) at (8,1) {};\node at (9.5,1) {$v(u_3)$};
\node[whitenode] (4) at (12,1) {};\node at (13.5,1) {$v(u_4)$};
\node at (16,0) {$\cdots$};
\node[whitenode] (5) at (18,1) {};\node at (19.5,1) {$v(u_5)$};
\node at (22,0) {$\cdots$};
\node[whitenode] (6) at (24,1) {};\node at (25.5,1) {$v(u_z)$};


\draw[dashed, thick] (0.5,4.3) rectangle (23.5,5.6) node[left = 0.2cm of 1'] {$W = W(\calT)$};


\node[blacknode] (U_1) at (14,-3) {};
\node[blacknode] (U_2) at (13,-3) {};
\node[blacknode] (U_3) at (12,-3) {};
\node[blacknode] (U_4) at (12,-6) {};
\node[blacknode] (U_5) at (13,-6) {};
\node[blacknode] (U_6) at (14,-6) {};
\node[whitenode] (U_0) at (15,-3) {};
\node[blacknode] (U_7) at (15,-6) {};

\node[blacknode] (U_8) at (13,-4) {};
\node[blacknode] (U_9) at (14,-4) {};
\node[blacknode] (U_10) at (13,-7) {};
\node[blacknode] (U_11) at (14,-7) {};
\node[whitenode] (U_12) at (16,-3) {};
\node[whitenode] (U_13) at (16,-6) {};
\node[whitenode] (U_14) at (15,-7) {};
\node[whitenode] (U_15) at (15,-4) {};

\draw[dashed, thick] (-0.5,0) rectangle (24.5,2) node[below right = 0.5cm and 3.5cm of 4] {base set (clique)};

\draw[-, thick, black!] (T1_1) -- (T1_2);
\draw[-, thick, black!] (T1_2) -- (T1_3);
\draw[-, thick, black!] (T1_3) -- (T1_4);
\draw[-, thick, black!] (T1_4) -- (T1_5);
\draw[-, thick, black!] (T1_5) -- (T1_6);
\draw[-, thick, black!] (T1_1) -- (T1_0);
\draw[-, thick, black!] (T1_6) -- (T1_7);

\draw[-, thick, black!] (T1_4) -- (T1_8);
\draw[-, thick, black!] (T1_8) -- (T1_9);
\draw[-, thick, black!] (T1_3) -- (T1_10);
\draw[-, thick, black!] (T1_10) -- (T1_11);
\draw[-, thick, black!] (T1_0) -- (T1_15);
\draw[-, thick, black!] (T1_7) -- (T1_14);
\draw[-, thick, black!] (T1_11) -- (T1_13);
\draw[-, thick, black!] (T1_9) -- (T1_12);

\draw[-, thick, black!] (T2_1) -- (T2_2);
\draw[-, thick, black!] (T2_2) -- (T2_3);
\draw[-, thick, black!] (T2_3) -- (T2_4);
\draw[-, thick, black!] (T2_4) -- (T2_5);
\draw[-, thick, black!] (T2_5) -- (T2_6);
\draw[-, thick, black!] (T2_1) -- (T2_0);
\draw[-, thick, black!] (T2_6) -- (T2_7);

\draw[-, thick, black!] (T2_4) -- (T2_8);
\draw[-, thick, black!] (T2_8) -- (T2_9);
\draw[-, thick, black!] (T2_3) -- (T2_10);
\draw[-, thick, black!] (T2_10) -- (T2_11);
\draw[-, thick, black!] (T2_0) -- (T2_15);
\draw[-, thick, black!] (T2_7) -- (T2_14);
\draw[-, thick, black!] (T2_11) -- (T2_13);
\draw[-, thick, black!] (T2_9) -- (T2_12);

\draw[-, thick, black!] (T3_1) -- (T3_2);
\draw[-, thick, black!] (T3_2) -- (T3_3);
\draw[-, thick, black!] (T3_3) -- (T3_4);
\draw[-, thick, black!] (T3_4) -- (T3_5);
\draw[-, thick, black!] (T3_5) -- (T3_6);
\draw[-, thick, black!] (T3_1) -- (T3_0);
\draw[-, thick, black!] (T3_6) -- (T3_7);

\draw[-, thick, black!] (T3_4) -- (T3_8);
\draw[-, thick, black!] (T3_8) -- (T3_9);
\draw[-, thick, black!] (T3_3) -- (T3_10);
\draw[-, thick, black!] (T3_10) -- (T3_11);
\draw[-, thick, black!] (T3_0) -- (T3_15);
\draw[-, thick, black!] (T3_7) -- (T3_14);
\draw[-, thick, black!] (T3_11) -- (T3_13);
\draw[-, thick, black!] (T3_9) -- (T3_12);

\draw[-, thick, black!] (T4_1) -- (T4_2);
\draw[-, thick, black!] (T4_2) -- (T4_3);
\draw[-, thick, black!] (T4_3) -- (T4_4);
\draw[-, thick, black!] (T4_4) -- (T4_5);
\draw[-, thick, black!] (T4_5) -- (T4_6);
\draw[-, thick, black!] (T4_1) -- (T4_0);
\draw[-, thick, black!] (T4_6) -- (T4_7);

\draw[-, thick, black!] (T4_4) -- (T4_8);
\draw[-, thick, black!] (T4_8) -- (T4_9);
\draw[-, thick, black!] (T4_3) -- (T4_10);
\draw[-, thick, black!] (T4_10) -- (T4_11);
\draw[-, thick, black!] (T4_0) -- (T4_15);
\draw[-, thick, black!] (T4_7) -- (T4_14);
\draw[-, thick, black!] (T4_11) -- (T4_13);
\draw[-, thick, black!] (T4_9) -- (T4_12);

\draw[-, thick, black!] (T5_1) -- (T5_2);
\draw[-, thick, black!] (T5_2) -- (T5_3);
\draw[-, thick, black!] (T5_3) -- (T5_4);
\draw[-, thick, black!] (T5_4) -- (T5_5);
\draw[-, thick, black!] (T5_5) -- (T5_6);
\draw[-, thick, black!] (T5_1) -- (T5_0);
\draw[-, thick, black!] (T5_6) -- (T5_7);

\draw[-, thick, black!] (T5_4) -- (T5_8);
\draw[-, thick, black!] (T5_8) -- (T5_9);
\draw[-, thick, black!] (T5_3) -- (T5_10);
\draw[-, thick, black!] (T5_10) -- (T5_11);
\draw[-, thick, black!] (T5_0) -- (T5_15);
\draw[-, thick, black!] (T5_7) -- (T5_14);
\draw[-, thick, black!] (T5_11) -- (T5_13);
\draw[-, thick, black!] (T5_9) -- (T5_12);

\draw[-, thick, black!] (T1_3) -- (1');
\draw[-, thick, black!] (T2_3) -- (2');
\draw[-, thick, black!] (T3_3) -- (3');
\draw[-, thick, black!] (T4_3) -- (4');
\draw[-, thick, black!] (T5_3) -- (5');

\draw[-, thick, black!] (U_1) -- (U_2);
\draw[-, thick, black!] (U_2) -- (U_3);
\draw[-, thick, black!] (U_3) -- (U_4);
\draw[-, thick, black!] (U_4) -- (U_5);
\draw[-, thick, black!] (U_5) -- (U_6);
\draw[-, thick, black!] (U_1) -- (U_0);
\draw[-, thick, black!] (U_6) -- (U_7);

\draw[-, thick, black!] (U_4) -- (U_10);
\draw[-, thick, black!] (U_8) -- (U_9);
\draw[-, thick, black!] (U_3) -- (U_8);
\draw[-, thick, black!] (U_10) -- (U_11);
\draw[-, thick, black!] (U_8) -- (U_15);
\draw[-, thick, black!] (U_10) -- (U_14);
\draw[-, thick, black!] (U_4) -- (U_13);
\draw[-, thick, black!] (U_0) -- (U_12);

\draw[-, thick, black!] (1) -- (U_3);
\draw[-, thick, black!] (2) -- (U_3);
\draw[-, thick, black!] (3) -- (U_3);
\draw[-, thick, black!] (4) -- (U_3);
\draw[-, thick, black!] (5) -- (U_3);
\draw[-, thick, black!] (6) -- (U_3);

\draw[-, thick, black!] (1) -- (1');
\draw[-, thick, black!] (1) -- (3');
\draw[-, thick, black!] (2) -- (1');
\draw[-, thick, black!] (2) -- (4');
\draw[-, thick, black!] (3) -- (2');
\draw[-, thick, black!] (3) -- (4');
\draw[-, thick, black!] (4) -- (1');
\draw[-, thick, black!] (4) -- (3');
\draw[-, thick, black!] (5) -- (2');
\draw[-, thick, black!] (5) -- (4');
\draw[-, thick, black!] (5) -- (5');
\draw[-, thick, black!] (6) -- (5');

\end{tikzpicture}
\caption{Instance of \minset{F}.} \label{fig:NP_F_all}
\end{subfigure}
\caption[]{Illustration of the graph $G_F$ constructed by Reduction~\ref{red:general}. Here $U = \{u_1, u_2, \ldots u_z\}$ and $\calT = \{T_1, T_2, \ldots , T_y\}$. The set of black vertices in Figure~\ref{fig:NP_C_all} is an example of a F-set of $G_F$.}
\label{fig:NP_F}
\end{figure}

We use the same vertex names as mentioned in Figure~\ref{fig:NP_F}. For an F-gadget (whose vertices are named as $b_j$ in Figure~\ref{fig:F-gadget}), the vertices of the graph $H_F(T)$ for some $T \in \calT$ or the graph $H_F(U)$ are denoted as $b_j(T)$ or $b_j(U)$, respectively. Moreover, for notational convenience, let $V_U = V(H_F(U))$ and let $V_T = V(H_F(T))$ for each $T \in \calT$. Next, we can show the following.

\begin{lemma} \label{lem:NP_F}
Let $A$ be an F-set of $G_F$. Then we have $\left |A \cap \Big ( V_U \cup \bigcup_{T \in \calT} V_T \Big ) \right | \ge 12|\calT|+11$.
\end{lemma}

\begin{proof}
Let $A$ be an F-set of $G_F$. For each $T \in \calT$, we must have $b_2(T), b_7(T) \in A$ in order for the pairs $b_1(T), b_{14}(T)$ and $b_8(T), b_{15}(T)$, respectively, to be I-separated (and hence, F-separated). Moreover, for each $T \in \calT$, we must have $b_3(T), b_6(T) \in A$ in order for the pairs $b_2(T), b_{14}(T)$ and $b_7(T), b_{15}(T)$, respectively, to be O-separated (and hence, F-separated). In addition, for each $T \in \calT$, we must have $b_9(T), b_{11}(T) \in A$ in order for the pairs $b_{10}(T), b_{13}(T)$ and $b_{12}(T), b_{16}(T)$, respectively, to be I-separated (and hence, F-separated). Also, among all the $(|\calT|+1)$ number of F-gadgets $(H_F,B_F)$ introduced by Reduction~\ref{red:general}, there can exist at most one support vertex and its adjacent leaf of one such $H_F$ which does not belong to $A$ (or else, there would be at least two leaves with empty closed neighborhoods in $A$ contradicting the fact that $A$ is an F-set of $G_F$). Thus, $A$ contains at least $12$ vertices from all $(|\calT|+1)$ number of $H_F$ introduced by Reduction~\ref{red:general} except possibly for one $H_F$ from which $A$ contains $11$ vertices. This proves the result.
\end{proof}

\begin{theorem} \label{thm:NP_F}
\minset{F} is NP-complete.
\end{theorem}

\begin{proof}
The problem clearly belongs to the class NP. Hence, to prove the result, we show that \minset{F} is also NP-hard. We do so by showing that an input $(U,\calT)$ to \test is a YES-instance of \test if and only if the input $(G_F,k)$ constructed in Reduction~\ref{red:general} is a YES-instance of \minset{F}.

($\implies$): Let $(U,\calT)$ be a YES-instance of \test. In other words, there exists a test collection $\calT'$ of $U$ such that $|\calT'| \le \ell$. We then construct a set $A \subset V(G_F)$ as follows. For each $T \in \calT'$, let $w(T) \in A$. Also, for all $T \in \calT$, let $b_i(T) \in A$ for all $i \in \{1, 2, \ldots, 12\}$. Similarly, let $b_i(U) \in A$ for all $i \in \{1, 2, \ldots, 12\} \setminus \{8\}$. This implies that $|A| = \ell + 12 |\calT|+11$. Since $\calT'$ is a test collection of $U$, by the construction of $G_F$, the set $A \cap W$ F-separates each pair of vertices in $Q$. Moreover, it can be verified that the vertices of $G_F$ have pairwise distinct open and pairwise distinct closed neighborhoods in $A$ (note that only $b_{15}(U)$ has an empty neighborhood in $A$). This makes $A$ an F-separating set of $G_F$ and hence, proves the necessary part of the result.

($\impliedby$): Let $(G_F,k)$ be a YES-instance of \minset{F}. In other words, there exists an F-set $A$ of $G$ such that $|A| \le k = \ell + 12 |\calT| + 11$. Now, using Lemma~\ref{lem:NP_F}, we must have $|A \cap W| \le \ell$. We now create a subcollection $\calT'$ of $\calT$ by letting $T \in \calT'$ if $w(T) \in A$. Then, we have $|\calT'| \le \ell$. Since $G_F[Q \cup \{b_5(U)\}]$ is a complete graph, the neighborhoods of the vertices $Q$ can differ only by their neighborhoods in $W$. Thus, pairs of vertices in $Q$ are F-separated
by the set $A \cap W$. Hence, by the construction of $G_F$, the collection $\calT'$ is a test collection of $U$. This proves the sufficiency part of the result and hence, proves the theorem.
\end{proof}

\subsection*{Locating sets}
\begin{figure}[!t]
\centering
\begin{subfigure}[t]{0.2\textwidth}
\centering
\begin{tikzpicture}[blacknode/.style={circle, draw=black!, fill=black!, thick, inner sep=2pt},
whitenode/.style={circle, draw=black!, fill=white!, thick, inner sep=2pt},
scale=0.4]

\hspace{-0.8cm}



\node[whitenode] (1) at (-1,8) {};\node at (-2,8) {$b_1$};
\node[whitenode] (2) at (3,8) {};\node at (4,8) {$b_2$};
\node[whitenode] (3) at (1,11) {};\node at (0,11) {$b_3$};
\node[whitenode] (4) at (1,14) {};\node at (0,14) {$b_4$};

\node () at (1,6) {$B_L = \{b_1, b_2\}$};

\node () at (0,2) {};

\draw[-, thick, black!] (1) -- (2);
\draw[-, thick, black!] (1) -- (3);
\draw[-, thick, black!] (2) -- (3);
\draw[-, thick, black!] (3) -- (4);


\end{tikzpicture}
\caption{L-gadget $(H_L, \{b_1,b2\})$.} \label{fig:L-gadget}
\end{subfigure}
\begin{subfigure}[t]{0.7\textwidth}
\centering
\begin{tikzpicture}[blacknode/.style={circle, draw=black!, fill=black!, thick, inner sep=1.5pt},
whitenode/.style={circle, draw=black!, fill=white!, thick, inner sep=1.5pt},
scale=0.4]
\scriptsize

\hspace{-0.8cm}

\node[blacknode] (T1_1) at (0,7) {};
\node[whitenode] (T1_2) at (2,7) {};
\node[blacknode] (T1_4) at (1,9) {};
\node[whitenode] (T1_5) at (1,11) {};

\node[blacknode] (T2_1) at (4,7) {};
\node[whitenode] (T2_2) at (6,7) {};
\node[blacknode] (T2_4) at (5,9) {};
\node[whitenode] (T2_5) at (5,11) {};

\node[blacknode] (T3_1) at (8,7) {};
\node[whitenode] (T3_2) at (10,7) {};
\node[blacknode] (T3_4) at (9,9) {};
\node[whitenode] (T3_5) at (9,11) {};

\node[blacknode] (T4_1) at (14,7) {};
\node[whitenode] (T4_2) at (16,7) {};
\node[blacknode] (T4_4) at (15,9) {};
\node[whitenode] (T4_5) at (15,11) {};

\node[blacknode] (T5_1) at (22,7) {};
\node[whitenode] (T5_2) at (24,7) {};
\node[blacknode] (T5_4) at (23,9) {};
\node[whitenode] (T5_5) at (23,11) {};

\node[blacknode] (1') at (1,5) {};\node at (2.5,5) {$w(T_1)$};0
\node[whitenode] (2') at (5,5) {};\node at (6.5,5) {$w(T_2)$};
\node[blacknode] (3') at (9,5) {};\node at (10.5,5) {$w(T_3)$};
\node at (13,5) {$\cdots$};
\node[blacknode] (4') at (15,5) {};\node at (16.5,5) {$w(T_i)$};
\node at (20,5) {$\cdots$};
\node[whitenode] (5') at (23,5) {};\node at (24.5,5) {$w(T_y)$};

\node[whitenode] (1) at (0,1) {};\node at (1.5,1) {$v^1(u_1)$};
\node[whitenode] (2) at (4,1) {};\node at (5.5,1) {$v^1(u_2)$};
\node[whitenode] (3) at (8,1) {};\node at (9.5,1) {$v^1(u_3)$};
\node[whitenode] (4) at (12,1) {};\node at (13.5,1) {$v^1(u_4)$};
\node at (16,0) {$\cdots$};
\node[whitenode] (5) at (18,1) {};\node at (19.5,1) {$v^1(u_5)$};
\node at (22,0) {$\cdots$};
\node[whitenode] (6) at (24,1) {};\node at (22.5,1) {$v^1(u_z)$};


\draw[dashed, thick] (0.5,4.2) rectangle (23.5,5.8) node[left = 0.2cm of 1'] {$W = W(\calT)$};


\draw[dashed, thick] (-0.5,0) rectangle (24.5,2); 

\node[whitenode] (1'') at (0,-4) {};\node at (1.6,-3.8) {$v^{\ell+1}(u_1)$};
\node[whitenode] (2'') at (4,-4) {};\node at (5.6,-3.8) {$v^{\ell+1}(u_2)$};
\node[whitenode] (3'') at (8,-4) {};\node at (9.6,-3.8) {$v^{\ell+1}(u_3)$};
\node[whitenode] (4'') at (12,-4) {};\node at (13.5,-3.8) {$v^{\ell+1}(u_4)$};
\node at (16,0) {$\cdots$};
\node[whitenode] (5'') at (18,-4) {};\node at (19.6,-3.8) {$v^{\ell+1}(u_5)$};
\node at (22,0) {$\cdots$};
\node[whitenode] (6'') at (24,-4) {};\node at (22.5,-3.7) {$v^{\ell+1}(u_z)$};

\node[blacknode] (p1) at (-2,0) {};\node at (-2.8,0) {$r^1_1$};
\node[whitenode] (p2) at (-2,2) {};\node at (-2.8,2) {$r^1_2$};
\node[blacknode] (p1') at (26,0) {};\node at (26.8,0) {$r^1_3$};
\node[whitenode] (p2') at (26,2) {};\node at (26.8,2) {$r^1_4$};


\draw[-, thick, black!] (p1) -- (-0.5,2);
\draw[-, thick, black!] (p1) -- (-0.5,0);
\draw[-, thick, black!] (p2) -- (-0.5,2);
\draw[-, thick, black!] (p2) -- (-0.5,0);

\draw[-, thick, black!] (p1') -- (24.5,2);
\draw[-, thick, black!] (p1') -- (24.5,0);
\draw[-, thick, black!] (p2') -- (24.5,2);
\draw[-, thick, black!] (p2') -- (24.5,0);

\draw[-, thick, black!] (p1) -- (p2);
\draw[-, thick, black!] (p1') -- (p2');


\node[blacknode] (q1) at (-2,-5) {};\node at (-2.8,-5) {$r^{\ell+1}_1$};
\node[whitenode] (q2) at (-2,-3) {};\node at (-2.8,-3) {$r^{\ell+1}_2$};
\node[blacknode] (q1') at (26,-5) {};\node at (26.8,-5) {$r^{\ell+1}_3$};
\node[whitenode] (q2') at (26,-3) {};\node at (26.8,-3) {$r^{\ell+1}_4$};


\draw[-, thick, black!] (q1) -- (-0.5,-3);
\draw[-, thick, black!] (q1) -- (-0.5,-5);
\draw[-, thick, black!] (q2) -- (-0.5,-3);
\draw[-, thick, black!] (q2) -- (-0.5,-5);

\draw[-, thick, black!] (q1') -- (24.5,-3);
\draw[-, thick, black!] (q1') -- (24.5,-5);
\draw[-, thick, black!] (q2') -- (24.5,-3);
\draw[-, thick, black!] (q2') -- (24.5,-5);

\draw[-, thick, black!] (q1) -- (q2);
\draw[-, thick, black!] (q1') -- (q2');


\node () at (-2,-1.25) {$\vdots$};
\node () at (0,-1.25) {$\vdots$};
\node () at (4,-1.25) {$\vdots$};
\node () at (8,-1.25) {$\vdots$};
\node () at (12,-1.25) {$\vdots$};
\node () at (18,-1.25) {$\vdots$};
\node () at (24,-1.25) {$\vdots$};
\node () at (26,-1.25) {$\vdots$};


\draw[dashed, thick] (-0.5,-5) rectangle (24.5,-3) node[below right = 2.5cm and 3cm of 4] {independent set};

\draw[dotted, thick, rounded corners] (-0.5,-4.6) rectangle (0.5,1.6);
\draw[dotted, thick, rounded corners] (3.5,-4.6) rectangle (4.5,1.6);
\draw[dotted, thick, rounded corners] (7.5,-4.6) rectangle (8.5,1.6);
\draw[dotted, thick, rounded corners] (11.5,-4.6) rectangle (12.5,1.6);
\draw[dotted, thick, rounded corners] (17.5,-4.6) rectangle (18.5,1.6);
\draw[dotted, thick, rounded corners] (23.5,-4.6) rectangle (24.5,1.6);

\node[blacknode] (U_1) at (11,-7) {};
\node[whitenode] (U_2) at (13,-7) {};
\node[whitenode] (U_4) at (12,-9) {};
\node[whitenode] (U_5) at (12,-11) {};

\draw[-, thick, black!] (T1_1) -- (T1_2);
\draw[-, thick, black!] (T1_1) -- (T1_4);
\draw[-, thick, black!] (T1_2) -- (T1_4);
\draw[-, thick, black!] (T1_4) -- (T1_5);

\draw[-, thick, black!] (T2_1) -- (T2_2);
\draw[-, thick, black!] (T2_1) -- (T2_4);
\draw[-, thick, black!] (T2_2) -- (T2_4);
\draw[-, thick, black!] (T2_4) -- (T2_5);

\draw[-, thick, black!] (T3_1) -- (T3_2);
\draw[-, thick, black!] (T3_1) -- (T3_4);
\draw[-, thick, black!] (T3_2) -- (T3_4);
\draw[-, thick, black!] (T3_4) -- (T3_5);

\draw[-, thick, black!] (T4_1) -- (T4_2);
\draw[-, thick, black!] (T4_1) -- (T4_4);
\draw[-, thick, black!] (T4_2) -- (T4_4);
\draw[-, thick, black!] (T4_4) -- (T4_5);

\draw[-, thick, black!] (T5_1) -- (T5_2);
\draw[-, thick, black!] (T5_1) -- (T5_4);
\draw[-, thick, black!] (T5_2) -- (T5_4);
\draw[-, thick, black!] (T5_4) -- (T5_5);

\draw[-, thick, black!] (U_1) -- (U_2);
\draw[-, thick, black!] (U_1) -- (U_4);
\draw[-, thick, black!] (U_2) -- (U_4);
\draw[-, thick, black!] (U_4) -- (U_5);

\draw[-, thick, black!] (T1_1) -- (1');
\draw[-, thick, black!] (T1_2) -- (1');
\draw[-, thick, black!] (T2_1) -- (2');
\draw[-, thick, black!] (T2_2) -- (2');
\draw[-, thick, black!] (T3_1) -- (3');
\draw[-, thick, black!] (T3_2) -- (3');
\draw[-, thick, black!] (T4_1) -- (4');
\draw[-, thick, black!] (T4_2) -- (4');
\draw[-, thick, black!] (T5_1) -- (5');
\draw[-, thick, black!] (T5_2) -- (5');

\draw[-, thick, black!] (1) -- (1');
\draw[-, thick, black!] (1) -- (3');
\draw[-, thick, black!] (2) -- (1');
\draw[-, thick, black!] (2) -- (4');
\draw[-, thick, black!] (3) -- (2');
\draw[-, thick, black!] (3) -- (4');
\draw[-, thick, black!] (4) -- (1');
\draw[-, thick, black!] (4) -- (3');
\draw[-, thick, black!] (5) -- (2');
\draw[-, thick, black!] (5) -- (4');
\draw[-, thick, black!] (5) -- (5');
\draw[-, thick, black!] (6) -- (5');

\draw[-, thick, black!] (1'') -- (U_1);
\draw[-, thick, black!] (1'') -- (U_2);
\draw[-, thick, black!] (2'') -- (U_1);
\draw[-, thick, black!] (2'') -- (U_2);
\draw[-, thick, black!] (3'') -- (U_1);
\draw[-, thick, black!] (3'') -- (U_2);
\draw[-, thick, black!] (4'') -- (U_1);
\draw[-, thick, black!] (4'') -- (U_2);
\draw[-, thick, black!] (5'') -- (U_1);
\draw[-, thick, black!] (5'') -- (U_2);
\draw[-, thick, black!] (6'') -- (U_1);
\draw[-, thick, black!] (6'') -- (U_2);

\end{tikzpicture}
\caption{Instance of \minset{L}.} \label{fig:NP_L_all}
\end{subfigure}
\caption[]{Illustration of the graph $G_L$ constructed by Reduction~\ref{red:general}. Here $U = \{u_1, u_2, \ldots u_z\}$ and $\calT = \{T_1, T_2, \ldots , T_y\}$. The vertices $v^1(u_i), \ldots , v^{\ell+1}(u_i)$ bordered with dotted vertical rectangles have the same neighborhood in $V(G_L) \setminus R$ although some edges have been omitted for a clearer figure. The set of black vertices in Figure~\ref{fig:NP_C_all} is an example of a L-set of $G_L$.}
\label{fig:NP_L}
\end{figure}
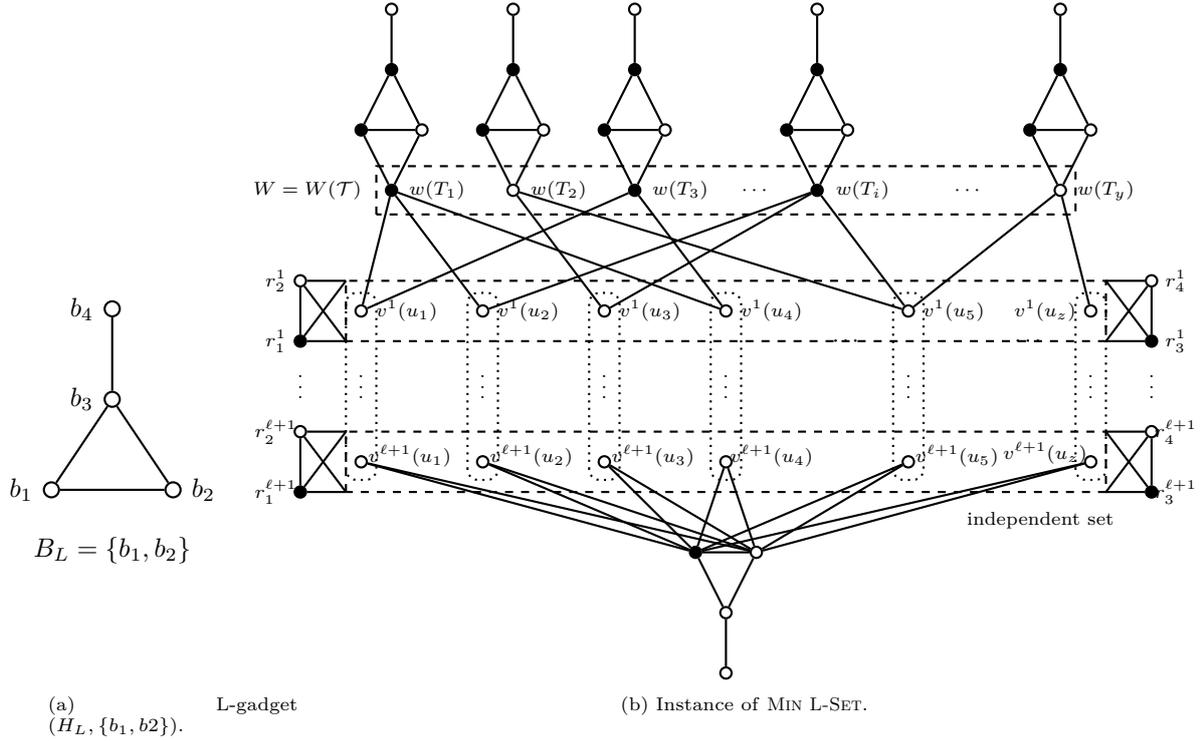

We next adapt Reduction~\ref{red:general} for the case $S = L$. In this case, in Step (1) of Reduction~\ref{red:general}, the reduction sets $r=r(L) = \ell + 1$ and so let $V(L,u) = \{v^1(u), v^2(u), \ldots , v^{\ell+1}(u)\}$ for each $u \in U$. Let $M = M(L) = \cup_{u \in U} V(S,u)$. Moreover, it sets $R = R(L) = \{r^i_1, r^i_2, r^i_3, r^i_4 : 1 \le i \le \ell+1\}$. The base set $Q = M \cup R$ is therefore of order $(\ell+1)(|U| + 4)$ and the reduction introduces edges so that the set $M$ is an independent set and each of the pairs $r^i_1, r^i_2$ and $r^i_3, r^i_4$ are closed twins each having neighborhood $\{v^i(u): u \in U\}$. The L-gadget in Step (3) of Reduction~\ref{red:general} is $(H_L, B_L)$ where $H_L$ is isomorphic to the graph shown in Figure~\ref{fig:L-gadget} and $B_L = \{b_1, b_2\}$ as shown in the figure. Finally, the reduction sets $p = p(L) = 2$ and $q = q(L) = 2\ell+3$ in Step (4) of Reduction~\ref{red:general}, that is, it sets $k = 3 \ell + 2 |\calT| + 3$.

We use the same vertex names as mentioned in Figure~\ref{fig:NP_L}. For an L-gadget (whose vertices are named as $b_j$ in Figure~\ref{fig:L-gadget}), the vertices of the graph $H_L(T)$ for some $T \in \calT$ or the graph $H_L(U)$ are denoted as $b_j(T)$ or $b_j(U)$, respectively. Moreover, for notational convenience, let $V_U = V(H_L(U))$ and let $V_T = V(H_L(T))$ for each $T \in \calT$. Next, we can show the following.

\begin{lemma} \label{lem:NP_L}
Let $A$ be an L-set of $G_L$. Then we have $\left |A \cap \Big ( R \cup V_U \cup \bigcup_{T \in \calT} V_T \Big ) \right | \ge 2|\calT|+2\ell+3$.
\end{lemma}

\begin{proof}
Let $A$ be an L-set of $G_L$. For each $T \in \calT$, since the pair $b_1(T), b_2(T)$ are closed twins in $G_L$, we must have either $b_1(T) \in A$ or $b_2(T) \in A$. Moreover, among all the $(|\calT|+1)$ number of L-gadgets $(H_L,B_L)$ introduced by Reduction~\ref{red:general}, there can exist at most one support vertex of one such $H_L$ which does not belong to $A$ (or else, there would be at least two leaves with empty neighborhoods in $A$ contradicting the fact that $A$ is an L-set of $G_L$). Thus, $A$ contains at least $2$ vertices from all $(|\calT|+1)$ number of $H_L$ introduced by Reduction~\ref{red:general} except possibly for one $H_L$ from which $A$ contains $1$ vertex. Finally, since the pairs $r^i_1, r^i_2$ and $r^i_3, r^i_4$ are twins in $G_L$ for all $i \in \{1, 2, \ldots , \ell+1\}$, we have $|\{r^i_1,r^i_2\} \cap A| \ge 1$ and $|\{r^i_3,r^i_4\} \cap A| \ge 1$. This proves the result.
\end{proof}

\begin{theorem} \label{thm:NP_L}
\minset{L} is NP-complete.
\end{theorem}

\begin{proof}
The problem clearly belongs to the class NP. Hence, to prove the result, we show that \minset{L} is also NP-hard. We do so by showing that an input $(U,\calT)$ to \test is a YES-instance of \test if and only if the input $(G_L,k)$ constructed in Reduction~\ref{red:general} is a YES-instance of \minset{L}.

($\implies$): Let $(U,\calT)$ be a YES-instance of \test. In other words, there exists a test collection $\calT'$ of $U$ such that $|\calT'| \le \ell$. We then construct a set $A \subset V(G_L)$ as follows. For each $T \in \calT'$, let $w(T) \in A$. Also, let $b_1(U) \in A$, let $r^i_1, r^i_3 \in A$ for all $i \in \{1, 2, \ldots , \ell+1\}$ and let $b_1(T), b_3(T) \in A$ for all $T \in \calT$. This implies that $|A| = 3\ell + 2 |\calT|+3$. Since $\calT'$ is a test collection of $U$, by the construction of $G_L$, the set $A \cap W$ L-separates each pair of vertices of the form $v^i(u), v^i(u')$, where $i \in \{1, 2, \ldots , \ell+1\}$ and $u, u' \in U$ with $u \ne u'$. On the other hand, pairs of the form $v^i(u), v^j(u')$ with $i \ne j$ are L-separated by the set $\{r^i_1, r^j_3\} \subset A$. Moreover, it can be verified that each vertex in $V(G_L) \setminus A$ has a unique open neighborhoods in $A$ (note that only $b_4(U)$ has an empty neighborhood in $A$). This makes $A$ an L-separating set of $G_L$ and hence, proves the necessary part of the result.

($\impliedby$): Let $(G_L,k)$ be a YES-instance of \minset{L}. In other words, there exists an L-set $A$ of $G$ such that $|A| \le k = 3\ell + 2 |\calT| + 3$. Moreover, we assume $A$ to be a minimum L-set of $G_L$. Now, using Lemma~\ref{lem:NP_L}, we must have $|A \cap (W \cup M)| \le \ell$ (recall that $M = M(L) = \cup_{u \in U} V(L,u)$). We now claim that $A \cap M = \emptyset$. To the contrary, let us assume that $v^i(u) \in A$ for some $i \in \{1, 2, \ldots , \ell+1\}$ and some $u \in U$. Notice that $v^i(u)$ may uniquely L-separate only itself from another vertex $v^i(u')$ for some $u' \in U$ with $u \ne u'$. If so, no vertex in $A \cap W$ L-separates $v^i(u)$ and $v^i(u')$. This implies that no vertex in $A \cap W$ L-separates the pairs $v^i(u), v^i(u')$ for all $i \in \{1, 2, \ldots , \ell+1\}$ (since the vertices $v^1(u), v^2(u), \ldots , v^{\ell+1}(u)$ are pairwise twins in $G_L$). Now, since $v^i(u)$ and $v^i(u')$ have the same neighborhood in $V(G_L) \setminus (A \cap W)$, we must have $|A \cap \{v^i(u), v^i(u')\}| \ge 1$ for all $i \in \{1, 2, \ldots , \ell+1\}$. However, using Lemma~\ref{lem:NP_L}, this implies that $|A| \ge 3\ell + 2|\calT| + 4$, contrary to our assumption that $|A| \le k = 3\ell + 2 |\calT| + 3$. Hence, the claim $A \cap M = \emptyset$ is true.

Thus, we have $A \cap (W \cup M) = A \cap W$ and so, $|A \cap W| \le \ell$. We now create a subcollection $\calT'$ of $\calT$ by letting $T \in \calT'$ if $w(T) \in A$. Then, we have $|\calT'| \le \ell$ as well. Since, any two vertices $v^1(u)$ and $v^1(u')$ for any $u,u' \in U$ are L-separated by some vertex in $W$, by the construction of $G_L$, the collection $\calT'$ is a test collection of $U$. This proves the sufficiency part of the result and hence, proves the theorem.
\end{proof}

\begin{proof}[Proof of Theorem~\ref{thm:NP_main}]
The proof follows by Theorems~\ref{thm_I-hard}, Corollary~\ref{cor:NP_O} and Theorems~\ref{thm:NP_F} and~\ref{thm:NP_L}.
\end{proof}

\section*{Proofs omitted in Section \ref{sec_rel}}

In this section, we established upper bounds on the X-number of a graph in terms of its S-number. 
We start with the proof of Theorem \ref{thm_S-2-SD-code}, stating that, for any $S \in \{L, O, I, F\}$ and an \sadmis ~graph $G$, we have $\sd(G) \leq \sset(G)+1$.

\begin{proof}[Proof of Theorem \ref{thm_S-2-SD-code}] 
Consider a (minimum) S-set $C \subseteq V$ of an \sadmis ~graph $G=(V,E)$. 
By definition of a dominating set, all vertices in $C$ are dominated. 
In $V \setminus C$, there is at most one vertex $v_0$ with $N[v_0] \cap C = \emptyset$:
\begin{itemize}
  \itemsep -3pt
  \item for $S \in \{I, F\}$, since $N[v] \cap C$ are unique sets for all $v \in V$ since $C$ is closed-separating; 
\item for $S \in \{L, O\}$, since $N[v] \cap C = N(v) \cap C$ holds for all $v \not\in C$ and $N(v) \cap C$ are unique sets for all $v \in V\setminus C$ (if $C$ is an L-set) respectively for all $v \in V$ (if $C$ is an O-set).
\end{itemize}
In all cases, either $C \cup \{v_0\}$ or $C$ (if no such vertex $v_0$ exists) is a dominating set (and still an S-set) of $G$, thus also an SD-code of $G$ and the assertion follows.
\end{proof}

We next wondered whether the same relation holds for the S-number of a graph and its STD-number.
Theorem \ref{thm_S-2-STD-code1} answered this question affirmatively for $S \in \{O, F\}$. 


\begin{proof}[Proof of Theorem \ref{thm_S-2-STD-code1}]
Consider a (minimum) S-set $C \subseteq V$ of an STD-admissible graph $G=(V,E)$ and $S \in \{O, F\}$. 
As $C$ is open-separating, there is at most one vertex $v_0$ with $N(v_0) \cap C = \emptyset$. 
Since $G$ is TD-admissible, hence has no isolated vertex, $v_0$ has (if it exists) a neighbor in $V \setminus C$, say $u_0$. 
We conclude that either $C \cup \{u_0\}$ or $C$ (if no vertex $v_0$ with $N(v_0) \cap C = \emptyset$ exists) is a total-dominating set (and still an S-set) of $G$, thus also an STD-code of $G$ and the assertion $\std(G) \leq \sset(G)+1$ follows.
\end{proof}

On the other hand, Theorem \ref{thm_S-2-STD-code2} showed that, for $S \in \{L, I\}$ and an \sadmis ~graph $G$, we have $\std(G) \leq 2 \sset(G)$. This is due to the following reasoning:

\begin{proof}[Proof of Theorem \ref{thm_S-2-STD-code2}]
Consider a (minimum) S-set $C \subseteq V$ of an STD-admissible graph $G=(V,E)$ and $S \in \{L, I\}$. 
In $V \setminus C$, there is at most one vertex $v_0$ with $N(v_0) \cap C = \emptyset$:
\begin{itemize}
  \itemsep -3pt
  \item for $S=L$, since $N(v) \cap C$ are unique sets for all $v \in V \setminus C$ since $C$ is locating; 
\item for $S=I$, since $N[v] \cap C = N(v) \cap C$ holds for all $v \not\in C$ and $N[v] \cap C$ are unique sets for all $v \in V$.
\end{itemize}
As $G$ is TD-admissible, hence has no isolated vertex, $v_0$ has (if it exists) at least one neighbor in $V \setminus C$, say $u_0$.  
Moreover, for every $v \in C$ with $N(v) \cap C = \emptyset$, there exists at least one neighbor in $V \setminus C$, say $u_v$ (again, as $G$ has no isolated vertex). Defining
$$
U = \{u_v \in V \setminus C : u_v \mbox{ is neighbor of $v \in C$ with $N(v) \cap C = \emptyset$} \},
$$
we see that $C \cup U \cup \{u_0\}$ is a TD-set of $G$. 
We next show that there is a TD-set $C^t \subseteq (C \cup U \cup \{u_0\})$ of cardinality $|C^t| \leq 2|C|$. 
The assertion clearly holds if
\begin{itemize}
  \itemsep -3pt
  \item there is $v \in C$ with $N(v) \cap C \neq \emptyset$ or there are $v,v' \in C$ having a common neighbor in $V \setminus C$ (as we then can choose $U$ with $|U| < |C|$);
\item there is no vertex $v_0 \in V \setminus C$ with $N(v_0) \cap C = \emptyset$ (as then adding $u_0$ is not needed).
\end{itemize}
Hence, assume that 
\begin{itemize}
  \itemsep -3pt
  \item we have $N(v) \cap C = \emptyset$ for all $v \in C$,
\item $U$ contains exactly one neighbor $u_v$ for each $v \in C$,
\item there is no vertex $u \in V \setminus C$ having two or more neighbors in $C$, and
\item there is a vertex $v_0 \in V \setminus C$ with $N(v_0) \cap C = \emptyset$.
\end{itemize}
Recall that $u_0$ is a neighbor of $v_0$ (as $v_0$  no isolated vertex). 
We have $u_0 \in V \setminus C$, and $u_0$ cannot be outside $U$: indeed, $u_0$ needs to have a neighbor in $C$ (as $v_0$ is the \emph{only} vertex with $N(v_0) \cap C = \emptyset$), but cannot have two or more neighbors in $C$ (by assumption). 
Hence, $u_0$ has exactly one neighbor in $C$, say $v'$. 
This implies that $u_0 = u_{v'}$ (as otherwise $u_0$ and $u_{v'}$ were not separated by $C$). 
Finally, $u_0 \in U$ implies that $C \cup U$ is the studied TD-set of $G$ in this case. 
Therefore, in all cases, there is a TD-set (still being an S-set) of cardinality $\leq 2|C|$.
\end{proof}

For $X \in \codes$, the X-problems have been studied from a unifying point of view, namely as reformulations in terms of covering problems in suitably constructed hypergraphs, for example, in \cite{ABLW_2018-DAM,ABLW_2022,ABW_2016,CW_ISCO2024,CW_Caldam2025}. 
Given a graph $G=(V,E)$ and an X-problem, we look for a hypergraph $\hyp_X(G) = (V,\hedge_X)$ so that $C \subseteq V$ is an X-code of $G$ if and only if $C$ is a \emph{cover} of $\hyp_X(G)$, that is, a subset $C \subseteq V$ satisfying $C \cap F \neq \emptyset$ for all $F \in \hedge_X$. Then the \emph{covering number} $\tau(\hyp_X(G))$, defined as the minimum cardinality of a cover of $\hyp_X(G)$, equals by construction the X-number $\x(G)$. The hypergraph $\hyp_X(G)$ is called the \emph{X-hypergraph} of $G$. 

It is a simple observation that encoding domination (respectively, total-domination) requires $\hedge_X$ to contain the closed (respectively, open) neighborhoods of all vertices of $G$. In order to encode the separation properties, that is, the fact that the intersection of an X-code with the neighborhood of each vertex is \emph{unique}, it was suggested in \cite{ABLW_2018-DAM,ABLW_2022} to use the symmetric differences of the neighborhoods. Here, given two sets $A$ and $B$, their \emph{symmetric difference} is defined by $A \bigtriangleup B = (A \setminus B) \cup (B \setminus A)$. 
In fact, it has been shown in \cite{ABLW_2018-DAM,ABLW_2022,CW_Caldam2025} that a subset $C \subseteq V$ of a graph $G=(V,E)$ is  
\begin{itemize}[leftmargin=12pt, itemsep=0pt]
  \itemsep -3pt
  \item locating if and only if $(N(u) \bigtriangleup N(v)) \cap C \neq \emptyset$ for all $uv \in E$ and $(N[u] \bigtriangleup N[v]) \cap C \neq \emptyset$ for all $uv \notin E$, where $u,v \in V$ are distinct, 
\item closed-separating if and only if $(N[u] \bigtriangleup N[v]) \cap C \neq \emptyset$ for all distinct $u,v \in V$,
\item open-separating if and only if $(N(u) \bigtriangleup N(v)) \cap C \neq \emptyset$ for all distinct $u,v \in V$, 
\item full-separating if and only if $(N[u] \bigtriangleup N[v]) \cap C \neq \emptyset$ for all $uv \in E$ and $(N(u) \bigtriangleup N(v)) \cap C \neq \emptyset$ for all $uv \notin E$, where $u,v \in V$ are distinct. 
\end{itemize}
We introduce analogously for all  $S \in \{L, O, I, F\}$ the \emph{separation-hypergraph} $\hyp_S(G) = (V,\hedge_S)$ so that $C \subseteq V$ is an S-set of $G$ if and only if $C$ is a cover of $\hyp_S(G)$. Then $\tau(\hyp_S(G))$ clearly equals by construction the S-number $\sset(G)$. 
From the previous findings on encoding the separation properties of X-codes for all $X \in \codes$, we immediately conclude that the 
hyperedge set $\hedge_S$ of $\hyp_S(G)$ consists of exactly 2 different subsets as shown in Table~\ref{tab_hypergraphs}, where \\[-5mm]
\begin{itemize}[leftmargin=12pt, itemsep=0pt]
\item $\bigtriangleup_a[G]$ (respectively, $\bigtriangleup_a(G)$) denotes the set of symmetric differences of closed (respectively, open) neighborhoods of all pairs of \emph{adjacent} vertices in $G$, 
\item $\bigtriangleup_n[G]$ (respectively, $\bigtriangleup_n(G)$) denotes the set of symmetric differences of closed (respectively, open) neighborhoods of all pairs of \emph{non-adjacent} vertices in $G$. \\[-5mm]
\end{itemize}
\begin{table}[ht]
\begin{center}
\begin{tabular}{ c || c | c | c | c  }
 S      & L                       & O                       & I                       & F \\ \hline 
 $\hedge_S$ & ~$\bigtriangleup_a(G)$~ & ~$\bigtriangleup_a(G)$~ & ~$\bigtriangleup_a[G]$~ & ~$\bigtriangleup_a[G]$~ \\
        & ~$\bigtriangleup_n[G]$~ & ~$\bigtriangleup_n(G)$~ & ~$\bigtriangleup_n[G]$~ & ~$\bigtriangleup_n(G)$ \\
\end{tabular}
\end{center}
\caption{
The S-hypergraphs $\hyp_S(G) = (V,\hedge_S)$ for the S-problems, listing column-wise the 2 needed subsets of hyperedges}
\label{tab_hypergraphs}
\end{table}

Consider the set $\mathcal{H}^*(V)$ of all hypergraphs $\hyp = (V,\mathcal{F})$ on the same vertex set $V$. 
We define a relation $\prec $ on $\mathcal{H}^*(V) \times \mathcal{H}^*(V)$ by $\hyp \prec \hyp'$ if, for every hyperedge $F'$ of $\hyp'$, there exists a hyperedge $F$ of $\hyp$ such that $F \subseteq F'$. In \cite{CW_Caldam2025}, it was shown that 
if $\hyp, \hyp' \in \mathcal{H}^*(V)$ with $\hyp \prec \hyp'$, then any cover of $\hyp$ is also a cover of $\hyp'$. In particular, we have $\tau(\hyp') \leq \tau(\hyp)$.
Combining this result with the hypergraph representations for the X-codes for all $X \in \codes$, the general relations of the X-numbers shown in Figure \ref{fig_relations} were obtained in \cite{CW_Caldam2025}. 
We immediately conclude for the S-numbers:

\begin{lemma}\label{lemA_S-relations}
For any accordingly admissible graph $G$, we have
$
\lset(G) \leq \oset(G),\iset(G) \leq \fset(G).
$
\end{lemma}


Moreover, it was noticed in e.g. \cite{ABLW_2018-DAM,ABLW_2022,ABW_2016,CW_ISCO2024,CW_Caldam2025} that $\hyp_X(G) = (V,\hedge_X)$ may contain redundant hyperedges. 
In fact, if there are two hyperedges $F, F' \in \hedge_X$ with $F \subseteq F'$, then $F \cap C \neq \emptyset$ also implies $F' \cap C \neq \emptyset$ for every $C \subseteq V$. 
Thus, $F'$ is \emph{redundant} as $(V,\hedge_X -\{F'\})$ suffices to encode the X-codes of $G$. 
Hence, only non-redundant hyperedges of $\hyp_X(G)$ need to be considered in order to determine $\tau(\hyp_X(G))$ and thus $\x(G)$. 
This clearly also applies to the here studied S-hypergraphs $\hyp_S(G) = (V,\hedge_S)$. 
Moreover, we call the hypergraph $\clu_S(G)$ obtained by removing all redundant hyperedges from $\hyp_S(G)$ the \emph{S-clutter} and note that clearly $\tau(\clu_S(G)) = \sset(G)$ holds.

In the sequel, it will turn out that several S-clutters are related to the hypergraph ${\mathcal{R}}_n^q=(V,{\cal F})$ called \emph{complete $q$-rose of order $n$}, where $V=\{1,\ldots,n\}$ and $\cal F$ contains all $q$-element subsets of $V$ for $2 \leq q < n$.  
From, e.g. \cite{ABLW_2018-DAM} we know that $\tau({\cal{R}}_n^q)=n-q+1$.
Note that, for $q=2$, ${\cal{R}}_n^q$ is in fact the clique 
$K_n$ and $\tau(K_n)=n-1$.

We are now prepared to prove Theorem \ref{thm_thin-spiders} providing the S-numbers of thin headless spiders $H_k=(Q \cup S, E)$ with $k \geq 4$, by
\[
\sset(H_k) = \left\{
\begin{array}{rl}
  k-1   & \mbox{if } S \in \{L,O\}  \\
  k+1   & \mbox{if } S = I,  \\
  2k-2  & \mbox{if } S = F.  \\
\end{array}
\right.
\]

\begin{proof}[Proof of Theorem \ref{thm_thin-spiders}]
Consider a thin headless spider $H_k=(Q \cup S, E)$ with $Q = \{q_1, \ldots, q_k\}$, $S = \{s_1, \ldots, s_k\}$ and $k \geq 4$.
In order to calculate $\sset(H_k)$ for $S \in \{L, O, I, F\}$, we construct the corresponding S-hypergraphs. 
To be able to build the symmetric differences, we first note that the closed and open neighborhoods are as follows: 
\begin{itemize}
  \itemsep -3pt
  \item $N[q_i] = Q \cup \{s_i\}$ and $N(q_i) = Q \setminus \{q_i\} \cup \{s_i\}$ for all $q_i \in Q$,
\item $N[s_i] = \{q_i,s_i\}$ and $N(s_i) = \{q_i\}$ for all $s_i \in S$.
\end{itemize}
This implies for the symmetric differences of closed or open neighborhoods of adjacent vertices, that 
\begin{itemize}
  \itemsep -3pt
  \item $N[s_i] \bigtriangleup N[q_i] = Q \setminus \{q_i\}$ and $N(s_i) \bigtriangleup N(q_i) = Q \cup \{s_i\}$ for $1 \leq i \leq k$,
\item $N[q_i] \bigtriangleup N[q_j] = \{s_i,s_j\}$ and $N(q_i) \bigtriangleup N(q_j) = \{s_i,s_j,q_i,q_j\}$ for $1 \leq i < j \leq k$,
\end{itemize}
and for the symmetric differences of closed or open neighborhoods of non-adjacent vertices, that 
\begin{itemize}
  \itemsep -3pt
  \item $N[s_i] \bigtriangleup N[q_j] = \{s_i,s_j\} \cup Q \setminus \{q_i\}$ and 
$N(s_i) \bigtriangleup N(q_j) = Q \setminus \{q_i,q_j\} \cup \{s_j\}$ for $i \neq j$,
\item $N[s_i] \bigtriangleup N[s_j] = \{s_i,s_j,q_i,q_j\}$ and 
$N(s_i) \bigtriangleup N(s_j) = \{q_i,q_j\}$ for $1 \leq i < j \leq k$.
\end{itemize}
Taking redundancies into account, we conclude the following.
\begin{itemize}
\item $\clu_O(H_k)$ is composed of only
  \begin{itemize}
  \item[] $N(s_i) \bigtriangleup N(s_j) = \{q_i,q_j\}$ for $1 \leq i < j \leq k$,
  \end{itemize}
  whereas all other hyperedges of $\hyp_O(H_k)$ are redundant. Hence $\clu_O(H_k)$ equals ${\mathcal{R}}_k^2$ induced by the vertices in $Q$ (whereas no vertex from $S$ is contained in any hyperedge) so that $\oset(H_k) = k-1$ follows.
\item $\clu_I(H_k)$ consists of 
  \begin{itemize}
  \itemsep -3pt
  \item[] $N[s_i] \bigtriangleup N[q_i] = Q \setminus \{q_i\}$ for $1 \leq i \leq k$ and 
  \item[] $N[q_i] \bigtriangleup N[q_j] = \{s_i,s_j\}$ for $1 \leq i < j \leq k$,
  \end{itemize}
  whereas all other hyperedges of $\hyp_I(H_k)$ are redundant. This shows that $\clu_I(H_k)$ has two components, ${\mathcal{R}}_k^{k-1}$ induced by the vertices in $Q$ and ${\mathcal{R}}_k^2$ induced by the vertices in $S$. Accordingly, we obtain $\oset(H_k) = \tau({\cal{R}}_k^{k-1}) + \tau({\cal{R}}_k^2) = 2 + k-1 = k+1$. 
\item $\clu_F(H_k)$ is composed of   
  \begin{itemize}
  \item[] $N(s_i) \bigtriangleup N(s_j) = \{q_i,q_j\}$ and $N[q_i] \bigtriangleup N[q_j] = \{s_i,s_j\}$ 
\end{itemize}
for $1 \leq i < j \leq k$, whereas all other hyperedges of $\hyp_F(H_k)$ are redundant. 
Thus $\clu_F(H_k)$ has two components equal to ${\mathcal{R}}_k^2$, one induced by the vertices in $Q$, the other induced by the vertices in $S$. Accordingly, we obtain $\fset(H_k) = 2 \tau({\cal{R}}_k^2) = 2(k-1) = 2k-2$.
\end{itemize}
Finally, from  
$\ld(G)  \leq \lset(G) +1$ and $\lset(G) \leq \oset(G)$ for all graphs $G$,  
$\ld(H_k) = k$ by \cite{ABLW_2022} and $\oset(H_k) = k-1$ from above, 
we conclude that also $\lset(H_k) = k-1$ follows.
\end{proof}

\section*{Proofs omitted in Section \ref{sec_coG}}

In this section, we addressed the question concerning the S-numbers of $G$ and its complement $\bar G$. 
We first prove the following:

\begin{lemma}\label{lemA_S_coG}
  For any graph $G$
  \begin{itemize}
    \itemsep -3pt
  \item we have $\hyp_L(G) = \hyp_L(\bar G)$;
  \item without closed twins, $\hyp_I(G) = \hyp_O(\bar G)$ holds;
  \item without open twins, we have $\hyp_O(G) = \hyp_I(\bar G)$;
  \item without twins, $\hyp_F(G) = \hyp_F(\bar G)$ follows.
  \end{itemize}
\end{lemma}

\begin{proof}
As we clearly have for any two distinct vertices $u,v$ of a graph $G$ that
$$
N_G(u) \bigtriangleup N_G(v) = N_{\bar G}[u] \bigtriangleup N_{\bar G}[v]
$$
holds, we immediately see that the hyperedge set 
\begin{itemize}
  \itemsep -3pt
  \item $\bigtriangleup_a(G)$ equals $\bigtriangleup_n[{\bar G}]$ and 
    \item $\bigtriangleup_n(G)$ equals $\bigtriangleup_a[{\bar G}]$.
\end{itemize}
Therefore, we conclude from Table~\ref{tab_hypergraphs} that
\begin{itemize}
  \itemsep -3pt
  \item $\bigtriangleup_a(G)$ and $\bigtriangleup_n[G]$ composing $\hyp_L(G)$ turn to 
        $\bigtriangleup_n[{\bar G}]$ and $\bigtriangleup_a({\bar G})$, i.e. constitute $\hyp_L(\bar G)$;
  \item $\bigtriangleup_a(G)$ and $\bigtriangleup_n(G)$ composing $\hyp_O(G)$ turn to 
        $\bigtriangleup_n[{\bar G}]$ and $\bigtriangleup_a[{\bar G}]$, i.e. yield $\hyp_I(\bar G)$;
  \item $\bigtriangleup_a[G]$ and $\bigtriangleup_n[G]$ composing $\hyp_I(G)$ turn to 
        $\bigtriangleup_n({\bar G})$ and $\bigtriangleup_a({\bar G})$, i.e. result in $\hyp_O(\bar G)$;
  \item $\bigtriangleup_a[G]$ and $\bigtriangleup_n(G)$ composing $\hyp_F(G)$ turn to 
        $\bigtriangleup_n({\bar G})$ and $\bigtriangleup_a[{\bar G}]$, i.e. yield $\hyp_F(\bar G)$.\\[-5mm]
\end{itemize}
\end{proof}
This implies particularly that 
  \begin{itemize}
    \itemsep -3pt
  \item an L-set of $G$ is also an L-set of $\bar G$;
  \item 
  an I-set of $G$ is an O-set of $\bar G$;
  \item 
  an O-set of $G$ is an I-set of $\bar G$;
  \item 
  an F-set of $G$ is also an F-set of $\bar G$;
\end{itemize}
which implies the assertion of Theorem \ref{thm_coS-numbers} that, for any graph $G$, we have
  \begin{itemize}
    \itemsep -3pt
  \item $\lset(G) = \lset(\bar G)$;
  \item $\iset(G) = \oset(\bar G)$;
  \item $\oset(G) = \iset(\bar G)$;
  \item $\fset(G) = \fset(\bar G)$.
  \end{itemize}

\end{document}